\documentclass[preprint]{elsarticle}
\usepackage{natbib}
\usepackage[all]{xy}
\usepackage{booktabs}

\usepackage{amssymb}
\usepackage{amsthm}
\usepackage{cite}
\usepackage{hyperref}
\usepackage{algorithm}
\usepackage{algpseudocode}
\usepackage[normalem]{ulem}
\usepackage{amsmath,mathtools}
\usepackage{amscd,enumitem}
\usepackage{caption}
\usepackage[justification=centering]{caption}
\usepackage{verbatim}
\usepackage{tabularx}
\usepackage{tabu}
\usepackage{multirow}
\usepackage{tikz-cd}
\usepackage{eurosym}
\usepackage{float}
\usepackage{color}
\usepackage{dcolumn}
\usepackage[mathscr]{eucal}
\usepackage[all]{xy}
\usepackage{bbm}
\usepackage{ytableau}
\usepackage[textheight=8.7in, textwidth=6.6in]{geometry}
\newtheorem*{conj*}{Conjecture}
\newtheorem{theorem}{Theorem}[section]

\theoremstyle{definition}
\newtheorem*{remark}{Remark}
\newtheorem*{tworemarks}{Two Remarks}
\newtheorem*{threeremarks}{Three Remarks}
\theoremstyle{plain}

\newtheorem{lemma}[theorem]{Lemma}

\newtheorem{proposition}[theorem]{Proposition}

\newtheorem{corollary}[theorem]{Corollary}
\newcommand{\legendre}[2]{\genfrac{(}{)}{}{}{#1}{#2}}
\newcommand{\Z}{\mathbb{Z}}

\newcommand{\Q}{\mathbb{Q}}
\newcommand{\R}{\mathbb{R}}

\newcommand{\oa}{\mathcal{O}}
\newcommand{\tr}{\mathrm{Tr}}

\numberwithin{equation}{section}

\newtheoremstyle{example}
  {\topsep}   
  {\topsep}   
  {\normalfont}  
  {0pt}       
  {\bfseries} 
  {.}         
  {5pt plus 1pt minus 1pt} 
  {}          
\theoremstyle{example}
\newtheorem*{example}{Example}

\usepackage{centernot}

\begin{document}
\begin{frontmatter}
\title{Arithmetic of Hecke $L$-Functions of Quadratic Extensions of Totally Real Fields}
\author[1]{Marie-H\'el\`ene Tom\'e}
\ead{mt392@duke.edu}
\affiliation[1]{organization={Department of Mathematics, Duke University},
            addressline={120 Science Dr},
            city={Durham},
            postcode={27708}, 
            state={NC},
            country={USA}}

\begin{abstract}
Deep work by Shintani in the 1970's describes Hecke $L$-functions associated to narrow ray class group characters of totally real fields $F$ in terms of what are now known as Shintani zeta functions. However, for $[F:\Q] = n \geq 3$, Shintani's method was ineffective due to its crucial dependence on abstract fundamental domains for the action of totally positive units of $F$ on $\R^n_+$, so-called \textit{Shintani sets}. These difficulties were recently resolved in independent work of Charollois, Dasgupta, and Greenberg and Diaz y Diaz and Friedman. For those narrow ray class group characters whose conductor is an inert rational prime in a totally real field $F$ with narrow class number $1$, we obtain a natural combinatorial description of these sets, allowing us to obtain a simple description of the associated Hecke $L$-functions. As a consequence, we generalize earlier work of Girstmair, Hirzebruch, and Zagier, that offer combinatorial class number formulas for imaginary quadratic fields, to real and imaginary quadratic extensions of totally real number fields $F$ with narrow class number $1$. For CM quadratic extensions of $F$, our work may be viewed as an effective affirmative answer to Hecke's Conjecture that the relative class number has an elementary arithmetic expression in terms of the relative discriminant.
\end{abstract}

\end{frontmatter}

\section{Introduction}\label{sec: intro}
    The $L$-functions associated to number fields are important tools in analytic number theory that provide information about the algebraic properties of these fields, such as class numbers, fundamental units, and regulators. A pioneering example is the work of Dirichlet, who in the case of quadratic fields $\Q(\sqrt{d})$ with fundamental discriminant $d$, gave the following class number formula (see \citep[Theorem~8.1.4]{koch})
        \begin{align} \label{dclf}
            h_{\Q(\sqrt{d})} = \begin{cases}
                    \frac{w_d\sqrt{|d|}}{2\pi} \cdot L(1, \chi_d) & \text{ if } d < 0 \\
                    \\[0.3pt]
                    \frac{\sqrt{d}}{\ln(\varepsilon_d)} \cdot L(1, \chi_d) & \text{ if } d > 0. 
            \end{cases}
        \end{align}
    Here $L(s, \chi_d)$ is the Dirichlet $L$-function associated to the Kronecker character $\legendre{d}{\cdot}$, $w_d$ is the number of roots of unity lying in $\Q(\sqrt{d})$, and $\varepsilon_d$ is the fundamental unit of $\Q(\sqrt{d})$.
    
    More recent class number formulas due to Girstmair \citep{girstmair}, Hirzebruch \citep{Hirzebruch}, and Zagier \citep{Zagier} elegantly simplify this formula in terms of familiar objects from elementary number theory; namely they make surprising connections with digit expansions and continued fractions. For primes $7 \leq p \equiv 3 \pmod{4}$ and $g$ a primitive root in $\mathbb{F}_p$, Girstmair showed that
        \begin{align} \label{g}
            h_{\Q(\sqrt{-p})} = \frac{1}{g+1}\sum_{k=1}^{p-1} (-1)^k x_k,
        \end{align}
    where $(x_1,x_2, \cdots, x_{p-1})$ are the digits of the periodic digit expansion of $1/p$ in base $g$. Under the additional assumption that $h_{\Q(\sqrt{p})} = 1$, Hirzebruch and Zagier proved
        \begin{align} \label{z}
            h_{\Q(\sqrt{-p})} = \frac{1}{3}\sum_{i=1}^{2t} (-1)^i a_i,
        \end{align}
    where $\sqrt{p}$ has continued fraction expansion $\sqrt{p} = [a_0, \overline{a_1a_2\cdots a_{2t}}]$. \footnote{It is a classical fact that these simple continued fractions have repeating digits of even period length.}
    
    These results follow from the fact that Dirichlet $L$-functions can be expressed as finite linear combinations of Hurwitz zeta functions, which satisfy functional equations relating their values at $s$ and $1-s$ and whose values at nonpositive integers can be expressed in terms of generalized Bernoulli numbers. Combining these facts, the formulas in (\ref{dclf}) become (see \citep[p. 234]{cohn})
        \begin{align} \label{nofear!}
            h_{\Q(\sqrt{d})} = \begin{cases}
                                    -\frac{w_d}{2|d|}\sum\limits_{r=1}^{|d|-1}\legendre{d}{r} r & \text{ if } d < 0 \\
                                    \\[0.3pt]
                                    -\frac{1}{2\ln(\varepsilon_d)} \sum\limits_{r=1}^{d-1}\legendre{d}{r}\ln\sin(\frac{\pi r}{d}) & \text{ if } d > 0.
                                \end{cases}
        \end{align}
    The alternating sums in (\ref{g}) and (\ref{z}) are combinatorial reformulations of the sums above.
    
    It is natural to ask whether class numbers can be expressed as finite alternating sums of combinatorial numbers beyond the realm of imaginary quadratic fields. We show that this is indeed the case for totally imaginary quadratic extensions of totally real fields $F$ with narrow class number $1$. 
    
    As the discussion above suggests, it is the combinatorial structure of $L$-functions themselves which underlies these formulas, and so we require a generalization of the theory of Dirichlet $L$-functions and their decomposition as sums of Hurwitz zeta functions. To this end, we consider certain Hecke $L$-functions of totally real fields of narrow class number $1$ and their combinatorial description in terms of generalizations of Hurwitz zeta functions. A significant step in this direction has already been obtained in the deep work of Shintani \citep{shintaniclnum, shintani}. We reformulate his work combinatorially in terms of number field invariants.
    
    We now recall the work of Shintani. Associated to a matrix $A \in \mathbb{M}_{n \times r}(\R_{>0})$ and a vector $\textbf{x} = (x_1, \cdots, x_r) \in \R_{\geq 0}^r$, Shintani defined \citep[$\S2$]{shintani} the \textit{Shintani zeta function}
        \begin{align} \label{Shintani zeta function}
        \zeta(s,A,\textbf{x}) \coloneqq \sum_{m_1, \ldots , m_r = 0}^{\infty}\prod_{i=1}^{n}
        \left(\sum_{j=1}^{r}a_{ij}(m_j + x_j)\right)^{-\textit{s}}, \quad \textrm{Re}(s) > 1
        \end{align}
    (in the notation of \citep{wei_lun_stark_units}), which coincides with the Hurwitz zeta function in the case $n=r=1$. He described (see Lemma 2 of \citep{shintani}) a decomposition of Hecke $L$-functions associated to totally real fields $F$ in terms of Shintani zeta functions evaluated at a finite set of algebraic points. His decomposition depends critically on the explicit description of the fundamental domain of the group action of the totally positive units of $\oa_F$ on $\R^n_+$, known as the \textit{Shintani set}. For totally real $F$ with $[F:\Q] \geq 3$, the construction of Shintani sets remained open. Therefore, Shintani's method offered a framework, albeit ineffective, for the decomposition of certain Hecke $L$-functions. Recent independent work by Charollois, Dasgupta, and Greenberg \citep{dasgupta} and Diaz y Diaz and Friedman \citep{ddf} filled the gap and provided explicit constructions of these Shintani sets.
    
    In view of this recent work, effective descriptions of Shintani's decomposition of Hecke $L$-functions can be obtained. Moreover, the relative class number formulas Shintani derived (see Theorem 2 of \citep{shintaniclnum}) for totally imaginary quadratic extensions of $F$ (i.e., where $[F:\Q] = n = 2$) become effective for $F$ of arbitrary degree over $\Q$. By carefully studying the combinatorial structure endowed by these Shintani sets, we obtain combinatorial descriptions of Shintani's decomposition of certain Hecke $L$-functions. This allows us to obtain reformulations of analogs of (\ref{nofear!}), in the spirit of Girstmair, Hirzebruch, and Zagier, for all quadratic extensions of $F$ with narrow class number $1$, where $[F:\Q]$ is arbitrary.
    
    To make this precise, we now turn to the problem of describing the Shintani sets. Recall that by Dirichlet's Unit Theorem, when $F$ is a totally real field of degree $n$ over $\Q$, the totally positive unit group $\oa_F^{\times, +}$ is a free abelian group of rank $n - 1$. Hence there exist $n-1$ totally positive units $\varepsilon_{1}, \cdots, \varepsilon_{n-1}$ such that $\oa_F^{\times, +} = \langle f_{\tau, 1}, \cdots, \varepsilon_{n-1} \rangle$. Following \citep{ddf, wei_lun_stark_units}, let the $n$ real embeddings of $F$ be given by $\sigma_1, \cdots, \sigma_n$ and let $\iota \coloneqq F \hookrightarrow \R^n$ be given by
        \begin{align} \label{i}
             x &\mapsto (\sigma_1(x), \cdots, \sigma_n(x)),\quad x \in F.
        \end{align}
    For any permutation $\tau \in S_{n-1}$, we define $f_{\tau, 1} \coloneqq 1$, and
        \begin{align} \label{gens}
            f_{\tau, j} \coloneqq \prod_{i=1}^{j-1} \varepsilon_{\tau(i)}, \quad 2 \leq j \leq n,
        \end{align}
    an associated matrix
        \begin{align} \label{mattau}
            A^{\tau} \coloneqq (\sigma_i(f_{\tau, j})) \in \mathbb{M}_n(\R_{>0}),
        \end{align}
   and a weight $w_\tau \in \{0, \pm 1\}$ (see Section \ref{sec: eff}). When $w_\tau$ is nonzero, the set of algebraic integers $\mathcal{B}_{F, \tau} \coloneqq \{f_{\tau, 1}, f_{\tau, 2}, \cdots, f_{\tau, n}\}$ forms a $\Q$-basis for $F$ and the set $\mathcal{B}_{\iota(F), \tau} \coloneqq \{\iota(f_{\tau,1}), \iota(f_{\tau, 2}), \cdots, \iota(f_{\tau, n})\}$ forms a basis for $\R^n$. Hence the lattice $\bigoplus_{i=1}^n \Z f_{\tau, i}$ is full rank.
   
   Let $e_n$ be the $n^{th}$ standard basis vector for $\R^n$, and denote by $(c_1, \cdots, c_n)$ the coefficients of $e_n$ in the basis $\mathcal{B}_{\iota(F), \tau}$, i.e., $e_n = \sum_{i=1}^n c_i\iota(f_{\tau, i})$. According to the sign of $c_i$, define $n$ intervals
        \begin{align} \label{ints}
                I_{\tau, i}\coloneqq
                \begin{cases}
                [0, 1) & \text{if $c_i >0$}\\
                (0, 1] & \text{otherwise},
                \end{cases} \quad 1 \leq i \leq n.
         \end{align}
    For any nonzero integral ideal $\mathfrak{f} \subset F$, the \textit{Shintani set} $R^{\tau}(\mathfrak{f})$ is defined by
        \begin{align} \label{Shintani Set}
            R^{\tau}(\mathfrak{f}) = \left\{z \in \mathfrak{f}^{-1} ~\bigg|~ z =\sum_{i=1}^n t_{z, \tau, i}f_{\tau, i}, \quad t_{z, \tau, i} \in I_{\tau, i} \right\}.
        \end{align}
    
    Now we turn to the problem of obtaining a combinatorial description of the Shintani sets associated to the ideals generated by inert rational primes. Throughout, we let $F \coloneqq \Q(\theta_F)$ for $\theta_F \in \oa_F$ be a totally real field with narrow class number $1$. Furthermore, we let $p \nmid [\oa_F:\Z[\theta_F]]$ be a rational prime which remains inert in $F$. Therefore, we have that $\oa_F/p\oa_F$ is isomorphic to $\mathbb{F}_{p^n}$ under the isomorphism $\varphi$ (see (\ref{phi})), and so we can fix $\rho$ such that $\mathbb{F}_{p^n}^{\times} = \langle \rho \rangle$. Let $h_\rho(x) \in \Z[x]$ be the minimal polynomial for $\rho$ whose reduction mod $p$ is a primitive polynomial in $\mathbb{F}_{p^n}$, say
        \begin{align} \label{minpoly}
            h_\rho(x) = x^n + p_{n-1}x^{n-1} + \cdots + p_0.
        \end{align}
    Using the coefficients of (\ref{minpoly}), we define a matrix $A_{F, \rho}(z)$ and a vector $\textbf{v}_{F, \rho}$, whose entries lie in the rational function field $\mathbb{Q}(z)$, by
        \begin{align} \label{mateq}
            A_{F, \rho}(z) \coloneqq \begin{pmatrix} 
            1 & 0 & 0 & \cdots & 0 & 0 & zp_0 \\
            -z & 1 & 0 & \cdots & 0 & 0 & zp_1 \\
             0 & -z & 1 & \cdots & 0 & 0 & zp_2 \\
            0 & 0 & -z & \cdots & 0 & 0 & zp_3 \\
            \vdots & \vdots & \vdots & \ddots & \vdots & \vdots & \vdots \\
            0 & 0 & 0 & \cdots & -z & 1 & zp_{n-2} \\
            0 & 0 & 0 & \cdots &0 & -z & 1 + zp_{n-1} 
            \end{pmatrix} \ \ \ \ \ \ \text{and} \ \ \ \ \ \  \textbf{v}_{F, \rho} \coloneqq \begin{pmatrix} -p_0 \\ -p_1 \\ -p_2 \\ -p_3 \\ \vdots \\ -p_{n-2} \\ -p_{n-1} \end{pmatrix}.
        \end{align}
    Since $\det(A_{F, \rho}(z)) = 1 + zp_{n-1} + z^2p_{n-2} + \cdots + z^np_0$ is a nonzero rational function (see Lemma \ref{nonsing}), by Cramer's Rule, there is a unique vector of rational functions with integer coefficients, say 
        \begin{align*} 
            \textbf{X}_{F, \rho} \coloneqq (X_{F, \rho, 1}(z), X_{F, \rho, 2}(z), \cdots, X_{F, \rho, n}(z)), 
        \end{align*}
    that satisfies $A_{F, \rho}(z)\textbf{X}_{F, \rho} = \textbf{v}_{F, \rho}.$ As power series, for $1 \leq i \leq n$, we have that 
        \begin{align} 
            X_{F, \rho, i}(z) = \sum_{m \geq 0} x_{F, \rho, i}(m)z^m = \sum_{m\geq 0} x_i(m)z^m,
        \end{align}
    where $x_i(0) = -p_i$. Note that for notational convenience, we drop the dependence on $F$ and $\rho$. For each permutation $\tau$ and each $1 \leq m < p^n$, the $n$-tuple $(x_{1}(m), \cdots, x_{n}(m))$ will be modified to obtain a finite set of $n$-tuples $\tilde{\textbf{x}}_\tau(i, m) \coloneqq (\tilde{x}_{\tau, 1}(i, m), \cdots, \tilde{x}_{\tau, n}(i, m))$, where $1 \leq i \leq \#\left(\oa_F \cap R^\tau(p\oa_F)\right)$. These $n$-tuples are the coefficients of elements of $R^\tau(p\oa_F)$ in the basis $\mathcal{B}_{F, \tau}$.

    Using the notation above, we explicitly describe Hecke $L$-functions for narrow ray class group characters with finite part of their conductor given by $p\oa_F$ in terms of the combinatorial data above and Shintani zeta functions.

    \begin{theorem}\label{thm: main thm}
        Assuming the notation and hypotheses above, we have that
            \begin{align*}
                    L(s, \chi_F) 
                    & = N(p\oa_F)^{-s} \sum_{\substack{\tau \in S_{n-1} \\ w_\tau \neq 0}} w_\tau \sum_{m=1}^{p^n-1} \mathrm{exp}\left(\frac{2\pi i k(n+m)}{d}\right) \sum_{i=1}^{\#\left(\oa_F \cap R^\tau(p\oa_F)\right)} \zeta\left(s, A^\tau, \tilde{\mathbf{x}}_{\tau}(i, m)\right),
            \end{align*}
        where $\chi_F$ is a narrow ray class group character with conductor $p\oa_F$, $N(p\oa_F)$ denotes the norm of the ideal $p\oa_F$, and $\mathrm{exp}\left((2\pi i k)/d\right)$, a primitive $d^{th}$ root of unity with $1 < d~|~(p^n-1)$, is the value of the finite part of the character $\chi_F$ on the equivalence class of $\varphi(\rho)$ in $(\oa_F/p\oa_F)^\times$.
    \end{theorem}

    \begin{tworemarks} \ \ \ 
        \newline
        \noindent
        (i) If $n =1$, then $F = \Q$ and every prime is inert in $F$. Therefore, Theorem \ref{thm: main thm} applies for every odd prime $p$, and gives the standard decomposition of $L(s, \chi_{d_p})$ as a linear combination of Hurwitz zeta functions where $d_p = (-1)^{\frac{p-1}{2}}p$.

        \smallskip
        \noindent
        (ii) If $F$ is a Galois extension of $\Q$ with cyclic Galois group, then by the Chebotarev Density Theorem, Theorem \ref{thm: main thm} applies to a positive density of primes. If $F/\Q$ is not Galois and there exists a rational prime $q$ which remains inert in $F$, then a positive density of primes remain inert in $F$, and hence Theorem \ref{thm: main thm} applies to a positive density of primes \citep[Proposition 8]{cueto-hernandez}.
    \end{tworemarks}

    Now we turn to the motivating problem of describing class numbers of imaginary quadratic extensions of totally real fields with narrow class number $1$.

    \begin{corollary} \label{cor: 1}
       Assume the notation and hypotheses above. Let $K = F(\sqrt{-p})$ and additionally assume that $p \equiv 3 \pmod{4}$. We have that
            \begin{align*}
                h_{K} &= \frac{1}{n} \cdot \frac{w_K}{[\oa_F^{\times}:\oa_F^{\times, +}][\oa_F^{\times, +}:N_{K/F}\oa_K^{\times}]} \sum_{\substack{\tau \in S_{n-1} \\ w_\tau \neq 0}} w_\tau & \\
                & \hspace{+2cm} \times  \left\{ \sum_{m=1}^{p^n-1}\sum_{i = 1}^{\#(\oa_F \cap R^\tau(p\oa_F))} (-1)^m \sum_{\substack{\tiny{(l_1,\dots,l_n)\in\mathbb{Z}_{\geq0}^{n}}\\\sum_{j=1}^{n}l_{j}=n}} \prod_{k=1}^n \frac{B_{l_k}(\tilde{x}_{\tau,k}(i, m))}{l_k!} \,\tr_{F/\Q}\bigg(\prod_{k=1}^n f_{\tau,k}^{l_k-1}\bigg)\right\},&
            \end{align*}
        where $w_K$ is the number of roots of unity lying in $K$, $N_{K/F}\oa_K^\times \coloneqq \{N_{K/F}(u)~|~u \in \oa_K^\times\}$, and $B_n(x)$ denotes the $n^{th}$ Bernoulli polynomial.
    \end{corollary}

    \begin{threeremarks} \ \ \ 
        \newline
        \noindent
        (i) When $n = 1$, the set $\{k/p : 1 \leq k < p\}$ considered by Girstmair is the Shintani set $R(p^{-1}\Z)$ generated by a primitive root $g$ of $\mathbb{F}_p^\times$ through the relation $x_1(m+1) = gx_1(m)$ for $1 \leq m < p$. This relation is encoded as the rational function $X_{F, g,1}(z) = gz/(1-gz).$

        \smallskip
        \noindent
        (ii) In the case where $F = \Q$ or $F$ is real quadratic, obtaining elementary formulas for the class number of a CM extension $K$ of $F$ is a classic problem appearing in the work of Dirichlet, Hecke, and Siegel \citep[p.130]{Siegel}. Hecke conjectured that the relative class number $h_K/h_F$ of CM extensions $K$ of totally real fields $F$ of arbitrary degree $n$ over $\Q$ also admits an elementary arithmetic description in terms of the relative discriminant of similar composition as (\ref{nofear!}) (see \citep[p. 2]{Hecke}). For a quadratic CM extension $K$ of a totally real field $F$ of this form, the above formula may be viewed as an effective affirmative answer to Hecke's Conjecture. 
        
        \smallskip
        \noindent
        (iii) Independently of this work, the case of Corollary (\ref{cor: 1}) when $n = 2$ was simultaneously obtained by Athaide, Cardwell, and Thompson in \citep{ECEpaper}.
    \end{threeremarks}
    
    Expressing $L(s, \chi_{F})$ in terms of evaluation of the Shintani zeta function over a union of Shintani sets sheds light on the combinatorics of the $L$-function and its higher order derivatives. The study of the Taylor expansion at $s = 0$ (or at $s=1$ via the functional equation) of the Hecke $L$-function $L(s, \chi_{F})$ associated to the narrow ray class group character $\chi_F$, is a fundamental problem in number theory. This expansion holds significant information for understanding arithmetic properties by using the data intrinsic to $F$. Hence, it is natural to investigate the higher order derivatives of Hecke $L$-functions associated to narrow ray class group characters, as they have the potential to encode essential arithmetic invariants.
    
    In particular, when $\chi_F$ is taken to be the narrow ray class group character associated to the field extension $K/F$ through class field theory, denoted $\chi_{K/F}$, the Taylor series expansion at $s=0$ of the associated $L$-function encodes arithmetic invariants of $K$. When $K/F$ is a totally imaginary extension, the constant term in the Taylor expansion gives the relative class number of $K/F$. When only one infinite place splits in $K$, the exponential of the first derivative of this $L$-function yields Stark units, which can be used to generate $K$. For $L$-functions associated to totally real quadratic extensions of $F$, the $n^{th}$ derivative of $L(s, \chi_{K/F})$ evaluated at $s=0$ gives the relative class number and also the relative regulator $R_K/R_F$ of the extension (for example, see \citep[Chapter 2]{BS}).
    
    \begin{corollary} \label{cor: 2}
        Assume the notation and hypotheses above. Let $K = F(\sqrt{p})$ and assume additionally that $p \equiv 1 \pmod{4}$. Then we have that
            \begin{align*}
                h_K\cdot\frac{R_K}{R_F} &= \frac{1}{n!} \sum_{k=0}^n (-1)^k{n \choose k} \left(\ln(N(p\oa_F))\right)^{n-k} \sum_{\substack{\tau \in S_{n-1}  \\ w_\tau \neq 0}} w_{\tau} \sum_{m = 1}^{p^n-1} (-1)^m \sum_{i=1}^{\#\left(\oa_F \cap R^\tau(p\oa_F)\right)} \zeta^{(k)} \left(0, A^\tau, \tilde{\mathbf{x}}_{\tau}(i, m)\right).
            \end{align*}
    \end{corollary}

    \begin{remark}
        Shintani also proved closed formulas for the first and second derivative of the Shintani zeta function evaluated at $s=0$ in terms of special functions. More precisely, in \citep[Proposition 1]{shintani}, he related the first derivative at $s=0$ to the logarithm of the Barnes multiple gamma function. However, not much is known about the values at $s=0$ of higher order derivatives.
    \end{remark}

    Now, we describe the organization of the paper. In Sections \ref{HLF} and \ref{SZF}, we recall basic results on Hecke $L$-functions and Shintani zeta functions. In Section \ref{sec: eff}, we recall the work of Charollois, Dasgupta, and Greenberg \citep{dasgupta} and Diaz y Diaz and Friedman \citep{ddf} which offers an effective construction of higher dimensional Shintani sets. Using these results, in Section \ref{sec: proof of thm 1.1}, we prove Theorem \ref{thm: main thm} and in Section \ref{sec: proof of corollaries}, we prove Corollaries \ref{cor: 1} and \ref{cor: 2} as a direct application of Theorem \ref{thm: main thm} to Hecke characters associated to relative quadratic extensions by class field theory. In Section \ref{section: examples}, we provide examples of the main application of Theorem \ref{thm: main thm} to computing class numbers of CM extensions of totally real fields $F$ of narrow class number 1 and $[F:\Q]$ arbitrary. To illustrate the explicit formulas obtained for the class number CM quadratic extensions $K/F$, we delve into detailed examples of Corollary \ref{cor: 1} for two cubic fields $F$.

\noindent 
\section{Hecke \textit{L}-functions and Shintani zeta functions}\label{sec: hlszf}
To prove Theorem \ref{thm: main thm}, we require some basic facts concerning Hecke $L$-functions, Shintani zeta functions, and Shintani sets which encode the interplay of these objects as dictated by the number fields in question.
    
\subsection{Hecke \textit{L}-Functions} \label{HLF}
    In this section, we follow \citep[Chapter VII $\S6$]{neukirch} (see \citep{Siegel} for a thorough discussion). Dirichlet $L$-functions encapsulate the behavior of rational primes in extensions of $\Q$. In the more general setting of number fields, the behavior of primes is captured by Hecke $L$-functions which naturally generalize Dirichlet $L$-functions. The simplest case of a Hecke $L$-function is the Dedekind zeta-function associated to a number field $F$, defined as
        \begin{align*}
            \zeta_F(s) = \sum_{\mathfrak{a} \subset \oa_F} N(\mathfrak{a})^{-s} = \prod_{\mathfrak{p}} (1-N(\mathfrak{p})^{-s})^{-1} \quad \mathrm{Re}(s) > 1,
        \end{align*}
    where $N(\mathfrak{a}) = |\oa_F/\mathfrak{a}|$, the sum ranges over all nonzero ideals of $\oa_F$, and the product is taken over prime ideals $\mathfrak{p}$ of $\oa_F$. Global information about $F$, such as the class number and the regulator, can be obtained from the analytic properties of the Dedekind zeta-function constructed from only local information about the primes of $\oa_F$. The analytic class number formula gives the following explicit formula for the residue of $\zeta_F(s)$ at $s=1$
        \begin{align} \label{gdclf}
            \mathop{\mathrm{Res}}_{s = 1} \zeta_F(s) = \frac{2^{r_1}(2\pi)^{r_2}h_FR_F}{\sqrt{|d_F|}{w_F}},
        \end{align}
     where $R_F$ is the regulator, $d_F$ is the absolute discriminant, $w_F$ is the number of roots of unity in $F$, and $r_1$ is the number of real embeddings of $F$ and $r_2$ is the number of pairs of complex embeddings of $F$. By the functional equation of the Dedekind zeta-function, (\ref{gdclf}) implies that the first nonzero term in the Taylor series expansion of $\zeta_F(s)$ around $s=0$ is given by (see \citep[p. 2]{Yoshida})
        \begin{align} \label{tayserdedzet}
            -\frac{h_FR_F}{w_F}s^{r_1+r_2-1}.
        \end{align}

    In the setting of $F = \Q$, for any positive integer $m$, a Dirichlet character mod $m$ is a character $\chi_m$ of the group $(\Z/m\Z)^{\times}$ extended to all $n \in \Z$ by letting $\chi(n) = 0$ when $(n,m) \neq 1$, and hence $\chi_m$ has period $m$ (see e.g., \citep[p. 27]{Davenport}). The Dirichlet $L$-function associated to $\chi_m$ is given by
        \begin{align*}
            L(s, \chi_m) = \sum_{n=1}^\infty \chi_m(n)n^{-s} = \prod_{p}(1-\chi_m(p)p^{-s})^{-1}, \quad \mathrm{Re}(s) > 1.
        \end{align*}
    
    Seeking a generalization of Dirichlet $L$-functions and their analytic properties, Hecke was led to the notion of what is now called a Hecke character. The $L$-functions associated to these characters can be shown to satisfy a functional equation. To make this precise, we require the following definitions. A \textit{modulus} $\mathfrak{m}$ is defined as the formal product of finite primes and real infinite primes of $F$. The modulus consists of
            \begin{align*}
                \mathfrak{m}_f = \prod_{\mathfrak{p} \nmid \infty} \mathfrak{p}^{\mathfrak{m}(\mathfrak{p})}, \quad \mathfrak{m}_{\infty} = \prod_{v | \infty} v^{\mathfrak{m}(v)},
            \end{align*}
    where $\mathfrak{m}(\mathfrak{p})\geq 0$ and $\mathfrak{m}(\mathfrak{p}) > 0$ for only finitely many finite primes of $F$, $\mathfrak{m}_f$ is an integral ideal of $\oa_F$, $\mathfrak{m}(v) \in \{0, 1\}$, and $\mathfrak{m}_{\infty}$ is a formal product of a subset of infinite real primes of $F$. A fractional ideal $\mathfrak{a} \subset I_F$ is said to be coprime to the modulus $\mathfrak{m}$ if no primes appearing in the decomposition of $\mathfrak{a}$ in $\oa_F$ appear in that of $\mathfrak{m}_f$. These ideals form a group, denoted $I_F^{\mathfrak{m}}$. A \textit{Hecke character} mod $\mathfrak{m}$ is a character $\chi: I_F^{\mathfrak{m}} \rightarrow \mathbb{S}^{1}$ such that there exists a pair of characters
        \begin{align*}
            \chi_f: (\oa_F/\mathfrak{m}_f)^{\times} \rightarrow \mathbb{S}^1 \quad \text{and} \quad \chi_{\infty}: \R^{\times} \rightarrow \mathbb{S}^1
        \end{align*}
    for which
            \begin{align*}
                \chi((\alpha)) = \chi_f(\alpha)\chi_{\infty}(\alpha),
            \end{align*}
    for every $\alpha \in \oa_F$ coprime to $\mathfrak{m}$. The character $\chi_f$ is a multiplicative function on $(\oa_F/\mathfrak{m}_f)^{\times}$, which extends to all invertible $\oa_F$-ideals by setting $\chi_f(\mathfrak{a}) = 0$ when $\mathfrak{a}$ is not coprime to the modulus $\mathfrak{m}$. The \textit{Hecke $L$-function} associated to a Hecke character $\chi$ is defined as
        \begin{align} \label{Lfgen}
            L(s, \chi) = \sum_{\mathfrak{a} \subset \oa_F} \chi(\mathfrak{a})N(\mathfrak{a})^{-s} = \prod_{\mathfrak{p}} \left(1-\chi(\mathfrak{p})N(\mathfrak{p})^{-s}\right)^{-1} \quad \mathrm{Re}(s) > 1,
        \end{align}
    where the sum is taken over all nonzero ideals of $\oa_F$ and the product is taken over all prime ideals $\mathfrak{p}$ of $\oa_F$.

    Throughout, we restrict to narrow ray class group characters defined as follows. Let
        \begin{align*}
            P_F^{\mathfrak{m}} \coloneqq \left\{(\alpha) \subset \oa_F~|~\alpha \equiv 1 \pmod{\mathfrak{m}} \text{ and } \sigma_v(\alpha) > 0 \text{ for all real infinite primes } v~|~\mathfrak{m}_{\infty}\right\},
        \end{align*}
   where $\sigma_v$ is the embedding associated to the infinite place $v$ and the congruence $\alpha \equiv 1 \pmod{\mathfrak{m}}$ means that $v_{\mathfrak{p}}(\alpha - 1) \geq v_{\mathfrak{p}}(\mathfrak{m}_f)$ for all primes $\mathfrak{p}$ dividing $\mathfrak{m}_f$. Then the \textit{ray class group mod} $\mathfrak{m}$ is the finite group $I_F^{\mathfrak{m}}/P_F^{\mathfrak{m}}$. If $\mathfrak{m}_{\infty}$ is the formal product of all the real infinite primes of $F$, the ray class group mod $\mathfrak{m}$ is called \textit{narrow}.

   Let $\mathfrak{m}$ be a modulus whose infinite part contains all real infinite places of $F$. If $\chi: I_F^{\mathfrak{m}} \rightarrow \mathbb{S}^1$ is a Hecke character for which
    \begin{align*}
                \chi\left((\alpha)\right) = \chi_f(\alpha)N\bigg(\bigg(\frac{\alpha}{|\alpha|}\bigg)^{\textbf{p}}\bigg),
    \end{align*}
    for all $\alpha \in \oa_F$, where $\chi_f$ is a character of $(\oa_F/\mathfrak{m}_f)^{\times}$ and $\textbf{p} \in \Z^{r_1+2r_2}$ is an admissible vector (see e.g., \citep[Chapter 7]{neukirch}), then $\chi$ is a \textit{narrow ray class group character} mod $\mathfrak{m}$. The \textit{conductor} of a narrow ray class group character $\chi$ mod $\mathfrak{m}$ is the smallest modulus $\mathfrak{f}$ dividing $\mathfrak{m}$ such that $\chi$ factors through the narrow ray class group $I_F^{\mathfrak{f}}/P_F^{\mathfrak{f}}$.
    
\subsection{Shintani Zeta Function} \label{SZF}
    The Barnes multiple zeta function is defined as
        \begin{align} \label{rprz}
            \zeta_r(s, \omega, x) \coloneqq \sum_{\Omega = m_1\omega_1 + \cdots + m_r\omega_r} (x + \Omega)^{-s},
        \end{align}
    where $\omega = (\omega_1, \cdots, \omega_r),$ $\omega_i > 0$ for all $1 \leq i \leq r$, $x > 0$, and $(m_1, \cdots, m_r)$ ranges over all $n$-tuples of non-negative integers (see e.g., \citep[p. 3]{Yoshida}). Barnes, using a method dating back to Riemann, proved that the Barnes multiple zeta function is holomorphic at $s = 0$ and has a meromorphic continuation to the whole complex plane with only simple poles at $s = 1, 2, \cdots, r$. The special value of the Barnes multiple zeta function at $s = 0$ has the following form (see e.g., \citep[p. 3]{Yoshida})
        \begin{align} \label{finform}
            \zeta_r\left(0, \omega, \sum_{k=1}^r \omega_k x_k\right) = (-1)^r \sum_{\substack{(l_1, \cdots, l_r)\in\Z_{\geq 0}^n \\ \sum_{i=1}^r l_r = r}} \omega_1^{l_1-1}\omega_2^{l_2-1}\cdots\omega_r^{l_r-1} \prod_{k=1}^r \frac{B_{l_k}(x_k)}{l_k!}.
        \end{align}
    
    Shintani generalized the definition of the Barnes multiple zeta function to a higher-dimensional zeta function taking as arguments an $n \times r$-matrix $A \in \mathbb{M}_{n\times r}(\R_{> 0})$ and a vector $\textbf{x} = (x_1, \cdots, x_r) \in \R_{\geq 0}^r$. Shintani \citep[$\S2$]{shintani} defined the Shintani zeta function 
         \begin{align}
        \zeta(s,A,\textbf{x}) \coloneqq \sum_{m_1, \ldots , m_r = 0}^{\infty}\prod_{i=1}^{n}
        \left(\sum_{j=1}^{r}a_{ij}(m_j + x_j)\right)^{-\textit{s}}, \quad \textrm{Re}(s) > 1.
        \end{align}
    The Shintani zeta function converges for $\mathrm{Re}(s) > r/n$, and has a meromorphic continuation to the whole complex plane. It is holomorphic except for possible poles at $s = 1, 2, \cdots, \lfloor r/n \rfloor$ and $s = t/n$ for integers $t$ satisfying $t \geq r$ and $n \nmid t$ (see \citep[p. 19]{Yoshida}). 
   
    When $n =1$, the Shintani zeta function associated to $\omega$, viewed as a $1 \times r$-matrix, and the scalar $x$ coincides with the Barnes multiple zeta-function. Therefore, for any row $A_i=(a_{i1}, a_{i2}, \cdots, a_{ir})$ of the matrix $A$, we have that (see e.g., \citep[p. 4]{Yoshida})
        \begin{align*}
            \zeta(s, A_i, \textbf{x}) = \zeta_r\left(s, (a_{i1}, a_{i2}, \cdots, a_{ir}), \sum_{j=1}^ra_{ij}x_j\right).
        \end{align*}
    Shintani \citep[$\S2$]{shintani} showed
        \begin{align*}
            \zeta(0, A, \textbf{x}) = \frac{1}{n}\sum_{i=1}^n \zeta(0, A_i, x) = \frac{1}{n}\sum_{i=1}^n\zeta_r\left(0, (a_{i1}, a_{i2}, \cdots, a_{ir}), \sum_{j=1}^ra_{ij}x_j\right).
        \end{align*}
    Using (\ref{finform}), he \citep{shintaniclnum} obtained a finite formula for the Shintani zeta function evaluated at $0$ (see e.g., \citep[Theorem 2.1]{Yoshida})
        \begin{align} \label{imp}
                \zeta(0, A, \textbf{x}) = \frac{(-1)^r}{n} \sum_{i=1}^n \sum_{\substack{\tiny{(l_1,\dots,l_n)\in\mathbb{Z}_{\geq0}^{n}}\\\sum_{j=1}^{r}l_{j}=r}} a_{i1}^{l_1-1}a_{i2}^{l_2-1}\cdots a_{ir}^{l_r -1} \frac{B_{l_1}(x_1)}{l_1!}\frac{B_{l_2}(x_2)}{l_2!}\cdots\frac{B_{l_r}(x_r)}{l_r!}.
            \end{align}
        
\subsection{Effective Shintani Set} \label{sec: eff}
     In what follows, we will have need for the following notation used in and results of Diaz y Diaz and Friedman \citep{ddf} and Barquero-Sanchez, Masri, and Tsai \citep{wei_lun_stark_units}. Recall that each permutation $\tau \in S_{n-1}$ determines a set $\{f_{\tau, 1}, \cdots, f_{\tau, n}\} \subset \oa_F$, defined as the products of the totally positive units $f_{\tau, 1}, \cdots, \varepsilon_{n-1}$ generating $\oa_F^{\times, +}$ as follows: $f_{\tau, 1} \coloneqq 1$ and for $2 \leq j \leq n$, $f_{\tau, j} \coloneqq \prod_{i=1}^{j-1} \varepsilon_{\tau(i)}$. For each $1\leq j \leq n$, the vector $\iota(f_{\tau, j}) \coloneqq (\sigma_1(f_{\tau, j}), \cdots, \sigma_n(f_{\tau, j})) \in \R^n$ is defined to be the vector of real embeddings of $f_{\tau, j}$. The $n \times n$-matrix with $\iota(f_{\tau, j})$ as its $j^{th}$ column is denoted by $A^\tau$. Using the matrix $A^\tau$ and the totally positive fundamental units $f_{\tau, 1}, \cdots, \varepsilon_n$, the weight $w_\tau$ associated to $\tau$ is defined as
        \begin{align*}
            w_\tau \coloneqq \frac{(-1)^{n-1} \text{sgn}(\tau) \cdot \text{sign}(\det(A^\tau))}{\text{sign}(\det(\log|\sigma_i(\varepsilon_j)|)_{1\leq i, j \leq n - 1})} \in \{0, \pm 1\}.
        \end{align*}
    For any nonzero integral ideal $\mathfrak{f} \subset F$, the \textit{Shintani set} $R^{\tau}(\mathfrak{f})$ is given by
        \begin{align*}
            R^{\tau}(\mathfrak{f}) = \left\{z \in \mathfrak{f}^{-1} ~\bigg|~ z =\sum_{i=1}^n t_{z, \tau, i}f_{\tau, i}, \quad t_{z, \tau, i} \in I_{\tau, i} \right\}.
        \end{align*}
    For notational convenience, let $\mathbf{t}_{z, \tau} = (t_{z, \tau, 1}, \cdots, t_{z, \tau, n})$. The following result of Diaz y Diaz and Friedman will be needed in the proof of Theorem \ref{thm: main thm}.
    \begin{corollary}[Corollary 3 of \citep{ddf}]
        Let $F$ be a totally real field of degree $n \geq 2$ with embeddings $\sigma_1, \cdots, \sigma_n$ and narrow class number $h_F^+ = 1$, let $\chi$ be a narrow ray class group character of $F$ and let the ideal $\mathfrak{f}$ be the finite part of its conductor. Then, for any set of generators $\varepsilon_1, \cdots, \varepsilon_{n-1}$ of the totally positive unit group of $F$, we have that
        \begin{align*}
            L(s, \chi) = N(\mathfrak{f})^{-s} \sum_{\substack{\tau \in S_{n-1} \\ w_\tau \neq 0}} w_\tau \sum_{z \in R^\tau(\mathfrak{f})} \chi(( z ) \mathfrak{f})\zeta(s, A^\tau, \mathbf{t}_{z,\tau}),
        \end{align*}
        where $(z)$ denotes the principal fractional ideal generated by $z \in F$, and 
        \begin{align*}
            \zeta(s,A^\tau,\textbf{t}_{z, \tau}) \coloneqq \sum_{m_1, \ldots , m_n = 0}^{\infty}\prod_{j=1}^{n}
        \left(\sum_{i=1}^n(t_{z, \tau, i} + m_i)\sigma_j(f_{\tau, i})\right)^{-\textit{s}}, \quad \mathrm{Re}(s) > 1
        \end{align*}
        is a Shintani zeta function.
    \end{corollary}
    The contribution from summing over each Shintani set $R^\tau(\mathfrak{f})$ depends on the permutation $\tau$ which determines $\{f_{\tau, 1}, \cdots, f_{\tau, n}\}$. The value of $w_\tau$ weights the contribution from the elements of $R^\tau(\mathfrak{f})$. The weight $w_\tau=0$ precisely when $A^\tau$ is singular, i.e., when the vectors $\{\iota(f_{\tau, i})~:~1\leq i \leq n\}$ do not form a basis of $\R^n$, and hence the totally positive units $\{f_{\tau, i}~:~1\leq i \leq n\}$ do not form a $\Q$-basis for $F$. When this is the case, the elements of $R^\tau(\mathfrak{f})$ do not contribute to the value of $L(s, \chi)$. In what follows, we assume $w_\tau \neq 0$. We denote the set $\{f_{\tau, 1}, \cdots, f_{\tau, n}\}$ giving a $\Q$-basis for $F$ by $\mathcal{B}_{F, \tau}$ and the set $\{\iota(f_{\tau, 1}), \cdots, \iota(f_{\tau, n})\}$ giving a basis for $\R^n$ by $\mathcal{B}_{\iota(F), \tau}$.
    
    Since $F$ is a $\Q$-vector space of rank $n$ with basis $\mathcal{B}_{F, \tau}$, every element of $F$ can be expressed as a unique $\Q$-linear combination of these basis elements. In particular, for any nonzero integral ideal $\mathfrak{f} \subset F$, any element $z$ lying in the fractional ideal $\mathfrak{f}^{-1}$ has a unique expression as a $\Q$-linear combination
        \begin{align*}
            z = \sum_{i=1}^n t_{z, \tau, i}f_{\tau, i}.
        \end{align*}
    The coordinates $t_{z, \tau, i}$ of $z \in \mathfrak{f}^{-1}$ in the basis $\mathcal{B}_{F, \tau}$ determine when $z$ lies in $R^\tau(\mathfrak{f})$. The set $R^\tau(\mathfrak{f})$ can be described as the set of all $z \in \mathfrak{f}^{-1}$ for which $\textbf{t}_{z, \tau}$ lies inside the bounded region of $\R^n$ determined by $I_{\tau, 1} \times \cdots \times I_{\tau, n}$, where these intervals are defined as in (\ref{ints}). Associated to each interval $I_{\tau, i}$ of the form $(0, 1]$ we define the following modified fractional part function
        \begin{align*}
            \{ x \}_{I_{\tau,i}} = 
            \begin{cases}
            \{ x \} & \text{if $x \not \in \Z$}\\
            1 & \text{if $x \in \Z$},
            \end{cases}
        \end{align*}
    otherwise, if $I_{\tau, i} = [0, 1)$, then $\{x\}_{I_{\tau, i}}$ is the usual fractional part function. It is readily seen that the coordinates of $z$ in the basis $\mathcal{B}_{F, \tau}$ and the coordinates of $\iota(z)$ in the basis $\mathcal{B}_{\iota(f_\tau)}$ coincide. The Shintani set $R^\tau(\mathfrak{f})$ may also be viewed as the subset of $\mathfrak{f}^{-1}$ whose image under $\iota$ falls within the fundamental parallelipiped $P^\tau$ for the full rank lattice $\bigoplus_{i=1}^n \Z \iota(f_{\tau, i}) \subset \R^n$. Hence the volume of $P^\tau$ is given by the determinant of $A^\tau$. Moreover, $R^\tau(\mathfrak{f})$ is finite because, as is the case for any fractional ideal of $F$, $\mathfrak{f}^{-1}$ determines a full rank lattice in $\R^n$. The quotient $G_\tau(\mathfrak{f}) \coloneqq \mathfrak{f}^{-1}/\bigoplus_{i=1}^n\Z f_{\tau, i}$ is a finite abelian group.

    Next we turn to the problem of describing the algebraic structure of the Shintani set. Throughout, we assume $\chi$ is a narrow ray class group character of $F$ with finite part of its conductor given by the principal ideal $\mathfrak{f} = ( \alpha )$.
    \begin{proposition}[Proposition 4.1 of\citep{wei_lun_stark_units}]
        $R^\tau(\mathfrak{f})$ is a complete set of coset representatives for the quotient group $G_\tau(\mathfrak{f})$ under the bijection
            \begin{eqnarray*}
                R^\tau(\mathfrak{f}) &\rightarrow& G_\tau(\mathfrak{f}) \\
                z &\mapsto& z + \bigoplus_{i=1}^n \Z f_{\tau, i}.
            \end{eqnarray*}
    \end{proposition} 
    \noindent By virtue of this bijection, one can define a binary operation $\oplus: R^\tau(\mathfrak{f}) \times R^\tau(\mathfrak{f}) \rightarrow R^\tau(\mathfrak{f})$ for any two elements $z_1, z_2 \in R^\tau(\mathfrak{f})$ using the group law of $G_\tau(\mathfrak{f})$. 
    \begin{proposition}[Proposition 4.3 of \citep{wei_lun_stark_units}]\label{binaryoperation} 
        The element $z_1 \oplus z_2$ is defined to be the unique coset representative of
            \begin{align*}
                z_1 + z_2 + \bigoplus_{i=1}^n \Z f_{\tau, i}
            \end{align*}
        lying in $R^\tau(\mathfrak{f})$. Under the group law $\oplus$, $R^\tau(\mathfrak{f})$ is a finite abelian group.
    \end{proposition} 
    \begin{remark} 
        The additive identity of $R^\tau(\mathfrak{f})$ is the element $1_{R^\tau(\mathfrak{f})} \in \bigoplus_{i=1}^n \Z f_{\tau, i}$ with coordinate vector $\textbf{t}_{z, \tau}$, where $t_{z, \tau, i} = 0$ or $1$ according to the definition of the interval $I_{\tau, i}$ (see (\ref{ints})). 
    \end{remark}
    Using the additive group structure of $R^\tau(\mathfrak{f})$, one can establish the following homomorphism of groups which will be needed in the proof of Theorem \ref{thm: main thm}.
    \begin{proposition}[Proposition 4.4 of \citep{wei_lun_stark_units}] \label{pi map}
        The map multiplication by $\alpha$ and reduction modulo $\mathfrak{f}$
            \begin{eqnarray*}
                \pi_{\alpha, \tau} \colon R^\tau(\mathfrak{f}) &\rightarrow \oa_F/\mathfrak{f}& \\
                z &\mapsto \alpha z + \mathfrak{f}&
            \end{eqnarray*} 
        is a surjective additive group homomorphism and for each coset $w + \mathfrak{f} \in \oa_F/\mathfrak{f}$,
                $$\#\pi_{\alpha, \tau}^{-1}(w + \mathfrak{f}) = \# \mathrm{ker}(\pi_{\alpha, \tau}).$$
    \end{proposition}
    \noindent By the First Isomorphism Theorem, $R^\tau(\mathfrak{f})/\mathrm{ker}(\pi_{\alpha, \tau}) \cong \oa_F/\mathfrak{f}$. Let 
        \begin{align} \label{bij}
            \pi^*_{\alpha, \tau}: R^\tau(\mathfrak{f})/\mathrm{ker}(\pi_{\alpha, \tau}) \xrightarrow{\sim} \oa_F/\mathfrak{f}
        \end{align}
    be the bijection on $R^\tau(\mathfrak{f})/\mathrm{ker}(\pi_{\alpha, \tau})$ induced by $\pi_{\alpha, \tau}$. The map $\pi^*_{\alpha, \tau}$ restricted to $(\oa_F/\mathfrak{f})^\times$ is a bijection between $R^\tau(\mathfrak{f})/\mathrm{ker}(\pi_{\alpha, \tau})-\{1_{R^\tau(\mathfrak{f})} + \mathrm{ker}(\pi_{\alpha, \tau})\}$ and $(\oa_F/\mathfrak{f})^\times$. For notational convenience, let 
        \begin{align*}
            R^\tau(\mathfrak{f})^\times \coloneqq R^\tau(\mathfrak{f})/\mathrm{ker}(\pi_{\alpha, \tau}) - \{1_{R^\tau(\mathfrak{f})} + \mathrm{\ker}(\pi_{\alpha, \tau})\},
        \end{align*}
    suggestive of the structure endowed by identifying $R^\tau(\mathfrak{f})^\times$ with $(\oa_F/\mathfrak{f})^\times$ under $\pi^*_{\alpha, \tau}$. Next we turn to an explicit characterization of the kernel of the map $\pi_{\alpha, \tau}$.
    \begin{lemma}[Lemma 4.6 \& Proposition 4.7 of\citep{wei_lun_stark_units}]\label{ker}
        We have that
            \begin{align*}
                \mathrm{ker}(\pi_{\alpha, \tau}) = \oa_F \cap R^\tau(\mathfrak{f}),
            \end{align*}
        and
            \begin{align*}
                \#\mathrm{ker}(\pi_{\alpha, \tau}) = \frac{\mathrm{vol}(P^\tau_F)}{\sqrt{d_F}} = \frac{\det (A^{\tau})}{\sqrt{d_F}},
            \end{align*}
        where $d_F$ is the discriminant of $F$.
    \end{lemma}

    Using the tools introduced in Section \ref{sec: hlszf} and the results of \citep{wei_lun_stark_units} we show the following lemma characterizing the relationship between narrow ray class group characters and the Shintani set for which we have need in the proof of Theorem \ref{thm: main thm}. We show that $\chi$ is invariant under translation by elements of $\mathrm{ker}(\pi_{\alpha, \tau})$.
    \begin{lemma} \label{trankerinv}
       Let $\chi$ be the narrow ray class group character which has the ideal $\mathfrak{f}$ as the finite part of its conductor. If $w \in \mathrm{ker}(\pi_{\alpha, \tau})$ then for any $z \in R^\tau(\mathfrak{f})$,
            \begin{align*}
                \chi((z) \mathfrak{f}) = \chi((z \oplus w) \mathfrak{f}).
            \end{align*}
    \end{lemma}
    \begin{proof}
        By Lemma \ref{ker}, $w \in \oa_F$ and hence $( w ) \mathfrak{f} \subset \mathfrak{f}$. By Proposition \ref{binaryoperation},
            \begin{align*}
                s \coloneqq z \oplus w - (z + w) \in \bigoplus_{i=1}^n \Z f_{\tau, i} \subset \oa_F,
            \end{align*}
        and so $w-s \in \oa_F$. Since $\mathfrak{f}$ is the finite part of the conductor of $\chi$,
            \begin{align*}
                \chi(( z \oplus w ) \mathfrak{f}) = \chi(( z + w - s) \mathfrak{f}) = \chi(( z ) \mathfrak{f} + ( w - s ) \mathfrak{f}) = \chi(( z ) \mathfrak{f}).
            \end{align*}
    \end{proof}

\section{Proof of Theorem 1.1}\label{sec: proof of thm 1.1}
    In this section, we will prove Theorem \ref{thm: main thm}. To this end, we first derive lemmas that will be required to prove Theorem 1.1 using the tools introduced in Section \ref{sec: hlszf}.
    \subsection{Some Lemmas}
    We specialize to the following setting. Let $K$ be a quadratic extension of a totally real field $F = \Q(\theta_F)$ of narrow class number $1$, where $\theta_F$ is chosen to be an algebraic integer. Further, suppose $p$ is a rational prime which remains inert in $F$ and $p \nmid [\oa_F: \Z[\theta_F]]$. Let $g(x) \in \Z[x]$ be the minimal polynomial of $\theta_F$ and denote by $\overline{g}(x) \in \mathbb{F}_p[x]$ its reduction mod $p$. Define the map
        \begin{align} \label{phi}
            \hspace{+2.5cm} \notag \varphi \colon \mathbb{F}_p[x]/(\overline{g}) &\rightarrow \oa_F/p\oa_F& \\
            x &\mapsto \theta_F + p\oa_F.&
        \end{align}
    
    Let $\chi_F$ be a narrow ray class group character with finite part of its conductor given by $p\oa_F$, where $p$ is a rational prime which remains inert in $F$. In this setting, the Shintani set $R^\tau(p\oa_F)$ has additional rich structure. Using the notation and tools developed in Section \ref{sec: hlszf}, we show that Shintani set $R^\tau(p\oa_F)$ is endowed with a cyclic structure mirroring the multiplicative structure of the multiplicative group of the finite field $\mathbb{F}_{p^n}$. This combinatorial structure underlies the proof of Theorem \ref{thm: main thm}.
    
    \begin{lemma} \label{bijec}
        Assume the notation and hypotheses above. Let $\varphi: \mathbb{F}_{p^n} \rightarrow \oa_F/p\oa_F$ be the map defined in (\ref{phi}) and let $\pi^*_{p, \tau}: R^\tau(p\oa_F)/\mathrm{ker}(\pi_{p, \tau}) \rightarrow \oa_F/p\oa_F$ be the bijection defined in (\ref{bij}). Then we have that the map
            \begin{align} \label{psi}
                \Psi \coloneqq (\pi^*_{p, \tau})^{-1} \circ \varphi \colon \mathbb{F}_{p^n} \rightarrow R^\tau(p\oa_F)/\mathrm{ker}(\pi_{p, \tau})
            \end{align}
        is a homomorphism of additive groups and $\mathbb{F}_{p^n}^\times$ and $R^\tau(p\oa_F)^\times$ are in bijective correspondence.
    \end{lemma}
    \begin{proof}
        Without loss of generality, fix $\tau \in S_{n-1}$ such that $w_\tau \neq 0$. By the discussion following Proposition \ref{pi map}, $R^{\tau}(p\oa_F)^\times$ and $(\oa_F/p\oa_F)^\times$ are in bijective correspondence under $\pi^*_{p, \tau}$.

        We have that $\Z[\theta_F]/p\Z[\theta_F] \cong \mathbb{F}_p[x]/(\overline{g})$. Since $p \nmid [\oa_F: \Z[\theta_F]]$ we may apply Dedekind's criterion to $p$ to obtain $\Z[\theta_F]/p\Z[\theta_F] \cong \oa_F/p\oa_F$ and hence since $p\oa_F$ is a prime ideal, $\varphi$ gives the following isomorphism of residue fields
            \begin{align*}
                \mathbb{F}_{p}[x]/(\overline{g}(x)) \cong \oa_F/p\oa_F
            \end{align*}
        via $x \mapsto \theta_F + p\oa_F$. Since $p$ remains inert, we have that the inertial degree $f_{p\oa_F | p} = [\oa_F:p\oa_F]$ equals $n$ and hence $\deg \overline{g}(x) = n$ so that
            \begin{align*}
                \mathbb{F}_{p^n} \cong \mathbb{F}_p(\theta_F) = \mathbb{F}_p[x]/(\overline{g}) \cong \oa_F/p\oa_F.
            \end{align*}
        Therefore, the composition $\Psi$ gives an additive group isomorphism from $R^\tau(p\oa_F)/\mathrm{ker}(\pi_{p, \tau})$ to $\mathbb{F}_{p^n}$. Since $1_{R^\tau(p\oa_F)} + \mathrm{ker}(\pi_{p, \tau})$ is the additive identity of $R^\tau(p\oa_F)/\mathrm{ker}(\pi_{p, \tau})$, the map \\ $\Psi|_{\mathbb{F}_{p^n}^\times} = (\pi^*_{p, \tau})^{-1} \circ \varphi|_{\mathbb{F}_{p^n}^{\times}}$ is a bijection between $\mathbb{F}_{p^n}^\times$ and $R^\tau(p\oa_F)^\times$.
    \end{proof}

    The multiplicative group of a finite field is cyclic, and so we may fix $\rho$ such that $\mathbb{F}_{p^n}^{\times} = \langle \rho \rangle$. Since $\mathbb{F}_{p^n}$ is the $(p^n-1)^{st}$ cyclotomic field over $\mathbb{F}_p$, the primitive elements $\rho$ of $\mathbb{F}_{p^n}$ are exactly the primitive $(p^n-1)^{st}$ roots of unity over $\mathbb{F}_p$. Hence the irreducible polynomial for a primitive element of $\mathbb{F}_{p^n}$ is one of the monic irreducible factors of the $(p^n-1)^{st}$ cyclotomic polynomial $\Phi_{p^n-1}(x)$ over $\mathbb{F}_p[x]$. Recall that we denoted by $h_\rho(x) = x^n + p_{n-1}x^{n-1} + \cdots + p_0 \in \Z[x]$ the minimal polynomial of $\rho$ over $\Q$ whose reduction mod $p$ is the minimal polynomial for $\rho$ over $\mathbb{F}_p$. Using the coefficients of $h_\rho$, we defined a matrix $A_{F, \rho}(z)$ lying in the rational function field $\mathbb{Q}(z)$ as
        \begin{align*}
            A_{F, \rho}(z) = \begin{pmatrix} 
            1 & 0 & 0 & \cdots & 0 & 0 & zp_0 \\
            -z & 1 & 0 & \cdots & 0 & 0 & zp_1 \\
             0 & -z & 1 & \cdots & 0 & 0 & zp_2 \\
            0 & 0 & -z & \cdots & 0 & 0 & zp_3 \\
            \vdots & \vdots & \vdots & \ddots & \vdots & \vdots & \vdots \\
            0 & 0 & 0 & \cdots & -z & 1 & zp_{n-2} \\
            0 & 0 & 0 & \cdots &0 & -z & 1 + zp_{n-1}
            \end{pmatrix}.
        \end{align*}
    \begin{lemma} \label{nonsing}
       The determinant of the $n \times n$ matrix $A_{F, \rho}(z)$ is given by
            \begin{align*} 
                \det(A_{F, \rho}(z)) = 1 + zp_{n-1} + z^2p_{n-2} + \cdots + z^np_0.
            \end{align*}
    \end{lemma}
    \begin{proof}
        We may write $A_{F, \rho}(z) = -z\left(C_{h_{\rho}} - z^{-1}I\right),$ where $C_{h_{\rho}}$ is the companion matrix for the polynomial $h_\rho(x)$. Therefore, since $\det(C_{h_{\rho}} -\lambda I) = (-1)^nh_\rho(\lambda)$ \citep[Theorem 1]{brand}, we have that
            \begin{align*}
                \det\left(-z\left(C_{h_{\rho}} - z^{-1}I\right)\right) = z^nh_\rho\left(z^{-1}\right)
                =1 + zp_{n-1} + z^2p_{n-2} + \cdots + z^np_0.&
            \end{align*}
    \end{proof}

    \subsection{Proof of Theorem 1.1}
        Define the map $\Phi: \mathbb{F}_{p^n} \rightarrow R^\tau(p\oa_F)/\mathrm{ker}(\pi_{p, \tau})$ as the composition of the following sequence of maps 
            \begin{align*}
            \Phi \coloneqq \mathbb{F}_{p^n} \xrightarrow{T_{\mathbb{F}_{p^n}, \mathcal{B}_{\rho}, \mathcal{B}_{\theta_F}}} \mathbb{F}_{p^n} \xrightarrow{\varphi} \oa_F/p\oa_F \xrightarrow{T_{F, \mathcal{B}_{\theta_F}, \mathcal{B}_{F, \tau}}} \oa_F/p\oa_F \xrightarrow{(\pi^*_{p, \tau})^{-1}} R^{\tau}(p\oa_F)/\mathrm{ker}(\pi_{p, \tau}),
            \end{align*}
        where $T_{\mathbb{F}_{p^n}, \mathcal{B}_{\rho}, \mathcal{B}_{\theta_F}}$ is an $\mathbb{F}_p$-vector space endomorphism taking the basis $\mathcal{B}_{\rho}$ to the basis $\mathcal{B}_{\theta_F}$ for $\mathbb{F}_{p^n}$ over $\mathbb{F}_p$ and $T_{F, \mathcal{B}_{\theta_F}, \mathcal{B}_{F, \tau}}$ is the vector space endomorphism taking the basis $\mathcal{B}_{\theta_F}$ to the basis $\mathcal{B}_{F, \tau}$ for $F$ over $\Q$. For notational convenience, we also denote the change of basis matrix in $\mathrm{GL}_n(\mathbb{F}_p)$ representing $T_{\mathbb{F}_{p^n}, \mathcal{B}_{\rho}, \mathcal{B}_{\theta_F}}$ by $T_{\mathbb{F}_{p^n}, \mathcal{B}_{\rho}, \mathcal{B}_{\theta_F}}$. Similarly, we use the same notation $T_{F, \mathcal{B}_{\theta_F}, \mathcal{B}_{F, \tau}}$ to denote the matrix in $\mathrm{GL}_n(\Q)$ representing the endomorphism $T_{F, \mathcal{B}_{\theta_F}, \mathcal{B}_{F, \tau}}$.
        
        Then we explicitly have that
            \begin{align*}
                \Phi(\rho^{n+m}) = \left(\sum_{i=1}^n \Bigg\{\frac{\big(T_{F, \mathcal{B}_{\theta_F}, \mathcal{B}_{F, \tau}}T_{\mathbb{F}_{p^n}, \mathcal{B}_{\rho}, \mathcal{B}_{\theta_F}}\overline{\textbf{x}}(m)\big)_{i}}{p}\Bigg\}_{I_{\tau, i}}f_{\tau, i}\right) + \mathrm{ker}(\pi_{p, \tau}),
            \end{align*}
        where $\overline{\textbf{x}}(m) = (\overline{x}_{1}(m), \overline{x}_{2}(m), \cdots, \overline{x}_{n}(m))$ is the minimal non-negative representative in $\mathbb{F}_{p^n}$ of the vector $(x_1(m), \cdots, x_n(m))$ holding the coefficients of $\rho^{n+m}$ in the basis $\mathcal{B}_{\rho}$ for $\mathbb{F}_{p^n}$. We will show that the $n$-tuple $(x_1(m), \cdots, x_n(m))$ is generated by $n$ unique rational functions $X_{F, \rho, 1}(z), \cdots, X_{F, \rho, n}(z)$ which can be explicitly solved for using the coefficients of the minimal polynomial $h_\rho(x)$ as described below. 

        Since $p$ remains inert in $F$ and $p \nmid [\oa_F:\Z[\theta_F]]$ by Dedekind's criterion, we have
            \begin{align*}
                \Z[\theta_F]/p\Z[\theta_F] \cong \mathbb{F}_p[x]/(\overline{g}(x)) = \mathbb{F}_p(\theta_F) \cong \oa_F/p\oa_F.
            \end{align*}
        Therefore $\mathcal{B}_{\theta_F}$ forms a power basis for $\mathbb{F}_{p^n}$ as a $\mathbb{F}_p$-vector space of dimension $n$. 

        For a primitive element $\rho$ of $\mathbb{F}_{p^n}^\times$, the set $\mathcal{B}_{\rho}$ is a power basis for $\mathbb{F}_{p^n}$ over $\mathbb{F}_p$. Hence there exists a change of basis matrix $T_{\mathbb{F}_{p^n}, \mathcal{B}_{\rho}, \mathcal{B}_{\theta_F}} \in \text{GL}_n(\mathbb{F}_p)$ taking an expression for $\rho^{n+m}$, $m > 0$ in the basis $\mathcal{B}_{\rho}$ to an expression of $\rho^{n+m}$ in the basis $\mathcal{B}_{\theta_F}$ by letting $T_{\mathbb{F}_{p^n}, \mathcal{B}_{\rho}, \mathcal{B}_{\theta_F}}$ act on the coefficient vector of $\rho^{n+m}$ for $m > 0$ in basis $\mathcal{B}_{\rho}$.
        
        Let $\varphi$ be defined as in (\ref{phi}) and let $\pi^*_{p, \tau}$ be the bijection defined in (\ref{bij}). Since $[F:\Q]=n$, we may view $F$  as a $\Q$-vector space of dimension $n$. Hence there exists a change of basis matrix $T_{F, \mathcal{B}_{\theta_F}, \mathcal{B}_{F, \tau}}$ taking the basis $\mathcal{B}_{\theta_F}$ to the basis $\mathcal{B}_{F, \tau}$.
        
        The powers of $\rho$ larger than $n-1$ can be written in the basis $\mathcal{B}_{\rho}$ using the relation given by $h_\rho$, namely,
            \begin{align*}
                \rho^n = -p_0 -p_1\rho- \cdots - p_{n-1}\rho^{n-1}
            \end{align*}
        so that
            \begin{align*}
                \rho^{n+1} = - p_{n-1}p_0 + \left(p_0 - p_{n-1}p_1\right)\rho + \cdots + \left(p_{n-2} - p_{n-1}p_{n-1}\right)\rho^{n-1}.
            \end{align*}
        Let $c_{i+1}(0) = -p_{i}$ for $i = 1, \cdots, n$ and for $m > 0$, assume the coefficients are properly chosen so that
            \begin{align} \label{recrel}
                \rho^{n+m} &= c_1(m) + c_2(m)\rho + \cdots + c_{n}(m)\rho^{n-1}.
            \end{align}
        Then we have that
            \begin{align} \label{relation}
                \notag \rho^{n+m+1} & = c_{1}(m)\rho + c_{2}(m)\rho^2 + \cdots + c_{n}(m)\rho^{n}& \\
                & = - c_{n}(m)p_0 + \left(c_{1}(m) - c_{n}(m)p_1\right)\rho + \cdots + \left(c_{n-1}(m) - c_{n}(m)p_{n-1}\right)\rho^{n-1},&
            \end{align}
        and hence by induction, for $1 \leq i \leq n$ and $m \geq 0$, the coefficients $c_i(m)$ satisfy (\ref{recrel}). Using (\ref{relation}) and the coefficients of $h_\rho$, we define the following matrix $A_{F, \rho}(z)$ and vector $\textbf{v}_{F, \rho}$, both lying in $\mathbb{Q}(z)$, by
        \begin{align*}
            A_{F, \rho}(z) \coloneqq \begin{pmatrix} 
            1 & 0 & 0 & \cdots & 0 & 0 & zp_0 \\
            -z & 1 & 0 & \cdots & 0 & 0 & zp_1 \\
            0 & -z & 1 & \cdots & 0 & 0 & zp_2 \\
            0 & 0 & -z & \cdots & 0 & 0 & zp_3 \\
            \vdots & \vdots & \vdots & \ddots & \vdots & \vdots & \vdots \\
            0 & 0 & 0 & \cdots & -z & 1 & zp_{n-2} \\
            0 & 0 & 0 & \cdots &0 & -z & 1 + zp_{n-1} 
            \end{pmatrix} \ \ \ \ \ \ \text{and} \ \ \ \ \ \  \textbf{v}_{F, \rho} \coloneqq \begin{pmatrix} 
            -p_0 \\
            -p_1 \\
            -p_2 \\
            -p_3 \\
            \vdots \\
            -p_{n-2} \\
            -p_{n-1} 
            \end{pmatrix}.
        \end{align*}
        Since $\det(A_{F, \rho}(z)) = 1 + zp_{n-1} + z^2p_{n-2} + \cdots + z^np_0$ is a nonzero rational function (see Lemma \ref{nonsing}), by Cramer's Rule, there exists a unique vector of rational functions with integer coefficients, say
            \begin{align*} 
                \textbf{X}_{F, \rho} \coloneqq (X_{F, \rho, 1}(z), X_{F, \rho, 2}(z), \cdots, X_{F, \rho, n}(z))
            \end{align*}
        that satisfies $A_{F, \rho}(z)\textbf{X}_{F, \rho} = \textbf{v}_{F, \rho}$. For notational convenience, we drop the dependence of the coefficients on $F$ and $\rho$, and let $x_{i}(m)$ be the $m^{th}$ coefficient of $X_{F, \rho, i}(z)$ as follows 
        \begin{align*} 
            X_{F, \rho, i}(z) = \sum_{m \geq 0} x_{F, \rho, i}(m)z^m = \sum_{m \geq 0}x_i(m) z^m.
        \end{align*}
        By Cramer's Rule, each of the rational functions $X_{F, \rho, 1}(z), \cdots, X_{F, \rho, n}(z)$ has denominator \\
        $1 + zp_{n-1} + \cdots + z^np_0$ of degree $n$. Note that $p_0 \neq 0$ since $h_\rho$ is irreducible. Setting $z = 0$, one verifies that for all $1 \leq i \leq n$, we have that $x_i(0)$ equals $c_i(0)$. Moreover, using this and the fact that $\mathbf{X}_{F, \rho}$ satisfies $A_{F, \rho}(z)\mathbf{X}_{F, \rho}=\mathbf{v}_{F, \rho}$, equating coefficients of $z^m$ one can verify that for all $1 \leq i \leq n$ and all $m \geq 0$, we have that $x_i(m)$ equals $c_i(m)$. 
        
        Define $\overline{\textbf{x}}(m) \coloneqq (\overline{x}_1(m), \cdots, \overline{x}_n(m))^T$ to be the minimal non-negative representative in $\mathbb{F}_p^n$ of the $n$-tuple of $m^{th}$ coefficients of each $X_{F, \rho, i}(z)$ viewed as a power series. We can perform a change of basis to express $\rho^{n+m}$ in the basis $\mathcal{B}_{\theta_F}$ for $\mathbb{F}_{p^n}$ corresponding to a modified vector $T_{\mathbb{F}_{p^n}, \mathcal{B}_{\rho}, \mathcal{B}_{\theta_F}} \overline{\textbf{x}}(m)$. Under the isomorphism $\varphi$, after applying the change of basis $T_{\mathbb{F}_{p^n}, \mathcal{B}_{\rho}, \mathcal{B}_{\theta_F}}$, we have that
            \begin{align*}
                \varphi(\rho^{n+m}) = (1, \theta_F, \cdots, \theta_F^{n-1}) \cdot T_{\mathbb{F}_{p^n}, \mathcal{B}_{\rho}, \mathcal{B}_{\theta_F}}\overline{\textbf{x}}(m)  + p\oa_F.
            \end{align*}
        Now applying the change of basis $T_{F, \mathcal{B}_{\theta_F}, \mathcal{B}_{F, \tau}}$, we obtain an expression for $\varphi(\rho^{n+m})$ in the basis $\mathcal{B}_{F, \tau}$. Namely, we have that
            \begin{align*}
                \varphi(\rho^{n+m}) = (f_{\tau, 1}, \cdots, f_{\tau, n})\cdot T_{F, \mathcal{B}_{\theta_F}, \mathcal{B}_{F, \tau}}T_{\mathbb{F}_{p^n}, \mathcal{B}_{\rho}, \mathcal{B}_{\theta_F}}\overline{\textbf{x}}(m) + p\oa_F.
            \end{align*}
        Under the mapping $(\pi^*_{p, \tau})^{-1}$, $\varphi(\rho^{n+m})$ is sent to the unique coset for $\varphi(\rho^{n+m})/p$ in $R^\tau(p\oa_F)/\mathrm{ker}(\pi_{p,\tau})$ given by
            \begin{align*}
                \left(\sum_{j=1}^{n} \bigg\{\frac{(T_{F, \mathcal{B}_{\theta_F}, \mathcal{B}_{F, \tau}}T_{\mathbb{F}_{p^n}, \mathcal{B}_{\rho}, \mathcal{B}_{\theta_F}}\overline{\textbf{x}}_{F, \rho}(m))_j}{p}\bigg\}_{I_{\tau, j}} f_{\tau, j}\right) + \mathrm{ker}(\pi_{p, \tau}).
            \end{align*}
        Since $\rho$ has order $p^n-1$ in $\mathbb{F}_{p^n}$, each element of the set $\{\rho^{n+m}~:~1\leq m\leq p^n-1\}$ is distinct. Since $\#(R^\tau(p\oa_F)/\mathrm{ker}(\pi_{p, \tau}))=\#(\mathbb{F}_{p^n})$ and $\Phi$ is an isomorphism, we have that
            \begin{align*}
               R^\tau(p\oa_F)/\mathrm{ker}(\pi_{p, \tau}) = \Phi(0_{\mathbb{F}_{p^n}}) \cup \mathrm{Im}_{\Phi}(\mathbb{F}_{p^n}^\times) = \{1_{R^\tau(p\oa_F)}\} \cup \{\Phi(\rho^{n+m}) : 1 \leq m \leq p^n-1\}.
            \end{align*}
        Therefore the set
            \begin{align*}
                C_\tau \coloneqq \{1_{R^\tau(p\oa_F)}\} \cup \left\{\vartheta(\rho^{n+m}) \coloneqq \sum_{i=1}^n\Bigg\{\frac{\big(T_{F, \mathcal{B}_{\theta_F}, \mathcal{B}_{F, \tau}}T_{\mathbb{F}_{p^n}, \mathcal{B}_{\rho}, \mathcal{B}_{\theta_F}}\overline{\textbf{x}}(m)\big)_{j}}{p}\Bigg\}_{I_{\tau, j}} f_{\tau, j} ~:~ 1 \leq m \leq p^n-1\right\}
            \end{align*}
        forms a complete set of coset representatives for the quotient group $R^\tau(p\oa_F)/\mathrm{ker}(\pi_{p, \tau})$.
        
        Put 
            \begin{align*} 
                \mathrm{ker}(\pi_{p, \tau}) = \left\{w(i)\coloneqq\sum_{j=1}^n w_j(i)f_{\tau, j}~:~1\leq i \leq \#\mathrm{ker}(\pi_{p, \tau})\right\}.
            \end{align*} 
        Then we have that
            \begin{align*}
                R^\tau(p\oa_F) - \{1_{R^\tau(p\oa_F)}\} = \{\vartheta(\rho^{n+m} ) \oplus w(i)~:~1 \leq i \leq \#\mathrm{ker}(\pi_{p, \tau}),  1 \leq m \leq p^n-1 \}.
            \end{align*}
        For $1 \leq j \leq n$, define the coefficients
            \begin{align*}
                \tilde{x}_{\tau, j}(i, m) \coloneqq \bigg\{ \bigg\{\frac{(T_{F, \mathcal{B}_{\theta_F}, \mathcal{B}_{F, \tau}}T_{\mathbb{F}_{p^n}, \mathcal{B}_{\rho}, \mathcal{B}_{\theta_F}}\overline{\textbf{x}}(m))_j}{p}\bigg\}_{I_{\tau, j}} + w_j(i) \bigg\}_{I_{\tau, j}}.
            \end{align*}
        Then by definition of the group law $\oplus$ we have that
            \begin{align*} 
                \vartheta(\rho^{n+m}) \oplus w(i) = \tilde{x}_{\tau, 1}(i, m)f_{\tau, 1} + \cdots + \tilde{x}_{\tau, n}(i, m)f_{\tau, n}.
            \end{align*}
        Hence the set of such $n$-tuples $(\tilde{x}_{\tau, 1}(i, m), \cdots, \tilde{x}_{\tau, n}(i, m))$, one for each $1 \leq i \leq \#\mathrm{ker}(\pi_{p, \tau})$ and $1 \leq m \leq p^n-1$, is in bijective correspondence with $R^\tau(p\oa_F) - \{1_{R^\tau(p\oa_F)}\}$ under the mapping
            \begin{align} \label{ssetc}
                (\tilde{x}_{\tau, 1}(i, m), \cdots, \tilde{x}_{\tau, n}(i, m)) \mapsto \tilde{x}_{\tau, 1}(i, m)f_{\tau, 1} + \cdots + \tilde{x}_{\tau, n}(i, m)f_{\tau, n}.
            \end{align}

        Since $\chi_F\left((1_{R^\tau(p\oa_F)})p\oa_F\right) = 0$ and by Lemma \ref{trankerinv}, the value of $\chi_F$ is invariant under translation by an element of the kernel, we may neglect the coset $1_{R^\tau(p\oa_F)} + \mathrm{ker}(\pi_{p, \tau})$. Moreover, by Lemma \ref{trankerinv},  the value of $\chi_F$ along a coset of $\mathrm{ker}(\pi_{p, \tau})$ is determined by the character value of a distinguished representative. Therefore, it suffices to consider the value of $\chi_F$ at each element of $C_\tau$.

        By Lemma \ref{ker}, $\vartheta(\rho^{n+m})$ differs from any element of the coset $\Phi(\rho^{n+m})$ by an algebraic integer. Therefore, 
            \begin{align*}
                \varphi(\rho^{n+m}) = (\pi^*_{p, \tau} \circ \Phi)(\rho^{n+m}) = p \cdot \vartheta(\rho^{n+m}) + p\oa_F.
            \end{align*}
        The restriction of $\varphi$ to $\mathbb{F}_{p^n}^\times$ gives a multiplicative group isomorphism between $\mathbb{F}_{p^n}^\times$ and $(\oa_F/p\oa_F)^{\times}$. Therefore, we have that
            \begin{align*}
                p \cdot \vartheta(\rho^{n+m}) + p\oa_F = \varphi(\rho^{n+m}) = \left(\varphi(\rho)\right)^{n+m} = \left(p \cdot \vartheta(\rho) + p\oa_F\right)^{n+m} = \left(p \cdot \vartheta(\rho)\right)^{n+m} + p\oa_F,&
            \end{align*}
        and hence
            \begin{align*}
                x \coloneqq p\cdot \vartheta(\rho^{n+m}) - p^{n+m} \cdot \vartheta(\rho)^{n+m} \in p\oa_F.
            \end{align*}
        Since $\chi_{F}$ is a multiplicative character of with finite part of its conductor given by $p\oa_F$,
            \begin{align*}
                \chi_{F}\left(( \vartheta(\rho^{n+m})) p\oa_F\right) = \chi_{F}\left((p\cdot\vartheta(\rho^{n+m}) - x) \oa_F\right) 
                = \left(\chi_F((\vartheta(\rho))p\oa_F)\right)^{n+m}.
            \end{align*}
            
        Since $\chi_F$ is a character of the narrow ray class group modulo $p\oa_F$, there exists a character \\ $\chi_f: (\oa_F/p\oa_F)^\times \rightarrow \mathbb{S}^1$ and a subset $S$ of $\{\sigma_1, \cdots, \sigma_n\}$ so that
            \begin{align*}
                \chi_F(x\oa_F) = \chi_f(x)\prod_{\sigma_i \in S} \mathrm{sgn}(\sigma_i(x)) 
            \end{align*}
        for any $x \in \oa_F$ (see e.g., \citep[p. 209]{shintani}). Since every element of the Shintani set is totally positive, for any $x \in R^\tau(p\oa_F)$, we obtain
            \begin{align*}
                \chi_F((px)\oa_F) = \chi_f(px),
            \end{align*}
        which yields
            \begin{align*} 
                \chi_F((p\cdot\vartheta(\rho))\oa_F) = \chi_f(p\cdot\vartheta(\rho)).
            \end{align*}
        Since $\varphi(\rho) = p\cdot\vartheta(\rho) + p\oa_F$, we have that $p\cdot \vartheta(\rho)$ lies in a nontrivial equivalence class mod $p\oa_F$ which generates $(\oa_F/p\oa_F)^\times$, and hence has order $p^n-1$. Moreover, since $\chi_f$ is a multiplicative homomorphism, a generator of $(\oa_F/p\oa_F)^\times$ is mapped to a primitive $d^{th}$ root of unity of order $d > 1$ dividing $p^n-1$.
        Therefore,  
             \begin{align*} 
                \chi_F((p\cdot \vartheta(\rho))) = \mathrm{exp}\left(\frac{2\pi i k}{d}\right),
            \end{align*}
       where $(k,d)=1$ and $d > 1$ is a divisor of $p^n-1$.

        By \citep[Corollary 3]{ddf},
                \begin{align*}
                    L(s, \chi_{F}) = N(p\oa_F)^{-s} \sum_{\substack{\tau \in S_{n-1} \\ w_\tau \neq 0}} w_\tau \sum_{z \in R^\tau(p\oa_F)} \chi_{F}(( z ) p\oa_F)\zeta(s, A^\tau, \mathbf{t}_{z, \tau}).
                \end{align*}
        Using (\ref{ssetc}) to run through the Shintani set and grouping by cosets of $\mathrm{ker}(\pi_{p, \tau})$ in $R^\tau(p\oa_F)$, we have that
                \begin{align*}
                    L(s, \chi_{F}) &= N(p\oa_F)^{-s}\sum_{\substack{\tau \in S_{n-1} \\ w_\tau \neq 0}} w_\tau \sum_{m=1}^{p^n-1}\sum_{i=1}^{\#\mathrm{ker}(\pi_{p, \tau})} \chi_{F}(( \vartheta(\rho^{n+m}) \oplus w_i ) p\oa_F)\zeta\left(s, A^\tau, \tilde{\mathbf{x}}_{\tau}(i, m)\right).
                \end{align*}
        By Lemma \ref{trankerinv},
                \begin{align*}
                    L(s, \chi_{F}) &= N(p\oa_F)^{-s} \sum_{\substack{\tau \in S_{n-1} \\ w_\tau \neq 0}} w_\tau \sum_{m=1}^{p^n-1}\sum_{i=1}^{\#\mathrm{ker}(\pi_{p, \tau})} \chi_{F}(( \vartheta(\rho^{n+m} ) p\oa_F)\zeta\left(s, A^\tau, \tilde{\mathbf{x}}_{\tau}(i, m)\right).&
                \end{align*}
        By Lemma \ref{ker}, we have that
                \begin{align*}
                    L(s, \chi_F) 
                    & = N(p\oa_F)^{-s} \sum_{\substack{\tau \in S_{n-1} \\ w_\tau \neq 0}} w_\tau \sum_{m=1}^{p^n-1} \chi_F((\vartheta(\rho)))^{n+m} \sum_{i=1}^{\#\left(\oa_F \cap R^\tau(p\oa_F)\right)} \zeta\left(s, A^\tau, \tilde{\mathbf{x}}_{\tau}(i, m)\right).
                \end{align*}
        Since $\chi_F((\varphi(\rho))) = \mathrm{exp}\left((2\pi i k)/d\right)$, we obtain
                \begin{align*}
                    L(s, \chi_F) 
                    & = N(p\oa_F)^{-s} \sum_{\substack{\tau \in S_{n-1} \\ w_\tau \neq 0}} w_\tau \sum_{m=1}^{p^n-1} \mathrm{exp}\left(\frac{2\pi i k(n+m)}{d}\right) \sum_{i=1}^{\#\left(\oa_F \cap R^\tau(p\oa_F)\right)} \zeta\left(s, A^\tau, \tilde{\mathbf{x}}_{\tau}(i, m)\right).
                \end{align*}

\section{Proof of Corollaries} \label{sec: proof of corollaries}
    In this section, we will prove Corollaries \ref{cor: 1} and \ref{cor: 2} using Theorem \ref{thm: main thm} and the tools introduced in Section \ref{sec: hlszf}. In order to prove Corollaries \ref{cor: 1} and \ref{cor: 2}, we require two lemmas.
         \begin{lemma} \label{reldiscim}
            Let $F$ be a totally real field with narrow class number $1$ and let $K = F(\sqrt{-p})$ where $p$ is a rational prime $p \equiv 3 \mod{4}$ which remains inert in $F$. Then $\mathfrak{D}_{K/F} = p\oa_F$.
        \end{lemma}
        \begin{proof}
            We show that $\{1, \frac{1 + \sqrt{-p}}{2}\}$ is an integral basis for $K$ over $F$. Since $h_F = 1$, $\oa_K$ is a free $\oa_F$-module of rank $2$ \citep[Corollary 3 on p. 386]{Narkiewicz}. Therefore, there exist $\omega_1, \omega_2 \in K$ such that $\oa_K = \omega_1\oa_F \oplus \omega_2\oa_F$. In particular, $\{\omega_1, \omega_2\}$ is an $F$-basis of $K$. Since $1 \in \oa_K$ and $\frac{1+\sqrt{-p}}{2} \in \oa_K$, there exist $m_{11}, m_{12}, m_{21}, m_{22} \in \oa_F$ such that
                \begin{align*}
                    m_{11} \omega_1 + m_{12}\omega_2 = 1 \\
                    m_{21}\omega_1 + m_{22}\omega_2  = \frac{1 + \sqrt{-p}}{2}.
                \end{align*}
            Namely,
                \begin{align*}
                    \begin{pmatrix} m_{11} & m_{12} \\
                    m_{21} & m_{22} \\
                    \end{pmatrix}
                    \begin{pmatrix} \omega_1 \\ \omega_2 \end{pmatrix} = \begin{pmatrix} 1 \\ \frac{1+\sqrt{-p}}{2} \end{pmatrix}
                \end{align*}
            where $M = (m_{ij}) \in \mathbb{M}_{2\times 2}(\oa_F)$.
            
            We now show that $M \in \text{GL}_2(\oa_F)$. If $M$ is invertible in $\oa_F$, then $M$ takes the integral basis of $\oa_F$ to another integral basis. By definition,
                \begin{align*}
                    \text{Disc}\bigg(\bigg\{1, \frac{1+\sqrt{-p}}{2}\bigg\}\bigg) = \det \begin{pmatrix}
                        \tr_{K/F}(1) & \tr_{K/F}\bigg(\frac{1+\sqrt{-p}}{2}\bigg) \\
                        \tr_{K/F}\bigg(\frac{1+\sqrt{-p}}{2}\bigg) & \tr_{K/F}\bigg(\bigg(\frac{1+\sqrt{-p}}{2}\bigg)^2\bigg) = -p\\
                    \end{pmatrix}.
                \end{align*}
            We also have that
                \begin{align*}
                    \text{Disc}\bigg(\bigg\{1, \frac{1+\sqrt{-p}}{2}\bigg\}\bigg)
                    & = \det \bigg(M \begin{pmatrix}
                        \tr_{K/F}(\omega_1^2) & \tr_{K/F}(\omega_1\omega_2) \\
                        \tr_{K/F}(\omega_2\omega_1) & \tr_{K/F}(\omega_2^2)
                        \end{pmatrix} M^T\bigg) & \\
                    & = \det (M)^2\cdot \mathrm{Disc}(\{\omega_1, \omega_2\}).&
                \end{align*}
            Hence we have 
                \begin{align*}
                    -p = \det(M)^2 \cdot \mathrm{Disc}(\{\omega_1, \omega_2\}).
                \end{align*}
            Since $p$ remains inert in $\oa_F$, $p$ cannot divide $\det(M)$. Since $F$ is a UFD, cancellation then shows $\det(M)$ must be a unit of $\oa_F$. Notice this also holds for primes which are unramified in $F$. Therefore, $M \in \text{GL}_2(\oa_F)$ and therefore $\{1, \frac{1 + \sqrt{-p}}{2}\}$ forms a $\oa_F$-integral basis of $K$. Thus we have that
                \begin{align*} 
                    \mathfrak{D}_{K/F} = \left(\det \begin{pmatrix}
                                1 && \frac{1+\sqrt{-p}}{2}\\
                                1 && \frac{1-\sqrt{-p}}{2}
                            \end{pmatrix}\right)^2 \oa_F =p\oa_F.
                \end{align*}
        \end{proof}

        \begin{lemma} \label{reldiscreal}
            Let $F$ a totally real field with $h_F^{+} = 1$ and let $K = F(\sqrt{p})$ where $p$ is a rational prime $p \equiv 1 \mod{4}$ which remains inert in $F$. Then $\mathfrak{D}_{K/F} = p\oa_F$.
        \end{lemma}
        \begin{proof}
            
        Applying the same argument as in Lemma (\ref{reldiscim}) to $\{1, \sqrt{p}\}$ shows $\{1, \sqrt{p}\}$ forms an $\oa_F$-integral basis of $K$, which yields
            \begin{align*} 
                \mathfrak{D}_{K/F} = \left(\det 
                \begin{pmatrix}
                    1 & \sqrt{p}\\
                    1 & \sqrt{p}
                \end{pmatrix}\right)^2\oa_F=p\oa_F.
        \end{align*}
        \end{proof}
    
    \subsection{Proof of Corollary 1.2} 
        When $\chi_{F}$ in Theorem \ref{thm: main thm} is taken to be the Hecke character $\chi_{K/F}$ associated to $K/F$ by class field theory, $\chi_{K/F}((z)p\oa_F) \in \{0, \pm 1\}$ for all $z \in F$. By the argument in the proof of Theorem \ref{thm: main thm}, $\chi_{K/F}((\vartheta(\rho))p\oa_F)$ is a primitive root of unity of order $d > 1$ dividing $p^n-1$. Therefore, $\chi_{K/F}((\vartheta(\rho))p\oa_F)$ equals $-1$, and hence we have that
            \begin{align} \label{charval--}
                \chi_{K/F}(( \vartheta(\rho^{n+m} )) p\oa_F) = (-1)^{n+m}.
            \end{align}
        By Theorem \ref{thm: main thm} applied to $\chi_{K/F}$, using (\ref{charval--}), we have that
            \begin{align*}
                L(0, \chi_{K/F}) = (-1)^n \sum_{\tau \in S_{n-1}} w_{\tau} \sum_{m=1}^{p^n-1} (-1)^m \sum_{i=1}^{\#\left(\oa_F \cap R^\tau(p\oa_F)\right)} \zeta\left(0, A^\tau, \tilde{\mathbf{x}}_{\tau}(i, m)\right).
            \end{align*}
        Now \citep[Theorem 2.1]{Yoshida} (see \ref{imp}) yields
            \begin{align*}
                \zeta(0, A^\tau, \tilde{\mathbf{x}}_\tau(i,m)) &= \frac{1}{n}\sum_{j=1}^n \zeta_n\bigg(0, (\sigma_j(f_{\tau, 1}), \sigma_j(f_{\tau, 2}), \cdots, \sigma_j(f_{\tau, n})), \sum_{k=1}^n \tilde{x}_{\tau, k}(i, m) \sigma_j(f_{\tau, k})\bigg)& \\
                & = \frac{(-1)^n}{n} \sum_{i=1}^n \sum_{\substack{\tiny{(l_1,\dots,l_n)\in\mathbb{Z}_{\geq0}^{n}}\\\sum_{j=1}^{n}l_{j}=n}} \sigma_i(f_{\tau, 1}^{l_1-1})\sigma_i(f_{\tau, 2}^{l_2-1})\cdots \sigma_i(f_{\tau, n}^{l_n -1}) \prod_{k=1}^n \frac{B_{l_k}(\tilde{x}_{\tau, k})}{l_k!}& \\
                & = \frac{(-1)^n}{n} \sum_{\substack{\tiny{(l_1,\dots,l_n)\in\mathbb{Z}_{\geq0}^{n}}\\\sum_{j=1}^{n}l_{j}=n}} \prod_{k=1}^n \frac{B_{l_k}(\tilde{x}_{\tau, k}(i, m))}{l_k!} \sum_{i=1}^n \prod_{k=1}^n \sigma_i(f_{\tau,k}^{l_k-1})& \\
                & = \frac{(-1)^n}{n}\sum_{\substack{\tiny{(l_1,\dots,l_n)\in\mathbb{Z}_{\geq0}^{n}}\\\sum_{j=1}^{n}l_{j}=n}} \prod_{k=1}^n \frac{B_{l_k}(\tilde{x}_{\tau, k}(i, m))}{l_k!} \, \tr_{F/\Q}\bigg(\prod_{k=1}^n f_{\tau, k}^{l_k-1}\bigg),
            \end{align*}
        where $\zeta_n(s, \omega, x)$ is the Barnes multiple zeta-function. Hence we have that
            \begin{align*}
                L(0, \chi_{K/F}) &= (-1)^n \sum_{\substack{\tau \in S_{n-1} \\ w_\tau \neq 0}} w_\tau \sum_{m=1}^{p^n-1} (-1)^m \sum_{i=1}^{\#\left(\oa_F \cap R^\tau(p\oa_F)\right)} \zeta\left(0, A^\tau, \tilde{\mathbf{x}}_\tau(i,m)\right)& \\
                & = (-1)^n \sum_{\substack{\tau \in S_{n-1} \\ w_\tau \neq 0}} w_\tau \sum_{m=1}^{p^n-1} \frac{(-1)^{m+n}}{n} \sum_{i=1}^{\#\left(\oa_F \cap R^\tau(p\oa_F)\right)} \sum_{\substack{\tiny{(l_1,\dots,l_n)\in\mathbb{Z}_{\geq0}^{n}}\\\sum_{j=1}^{n}l_{j}=n}} \prod_{k=1}^n \frac{B_{l_k}(\tilde{x}_{\tau, k}(i, m))}{l_k!} \tr_{F/\Q}\left(\prod_{k=1}^n f_{\tau,k}^{l_k-1}\right)& \\
                & = \frac{1}{n} \sum_{\substack{\tau \in S_{n-1} \\ w_\tau \neq 0}} w_\tau \sum_{m=1}^{p^n-1} (-1)^{m} \sum_{i=1}^{\#\left(\oa_F \cap R^\tau(p\oa_F)\right)} \sum_{\substack{\tiny{(l_1,\dots,l_n)\in\mathbb{Z}_{\geq0}^{n}}\\\sum_{j=1}^{n}l_{j}=n}} \prod_{k=1}^n \frac{B_{l_k}(\tilde{x}_{\tau, k}(i, m))}{l_k!} \tr_{F/\Q}\left(\prod_{k=1}^n f_{\tau,k}^{l_k-1}\right).&
            \end{align*}
            
        Next, we consider the factorization
        \begin{align*}
            \zeta_K(s) = \zeta_F(s)L(s, \chi_{K/F}).
        \end{align*}
        For the Dedekind zeta function associated to any number field $K$, the Taylor series expansion at $s = 0$ is given by
            \begin{align*}
                \zeta_K(s) = \frac{-h_KR_K}{w_K} s^{r_1 + r_2 - 1} + O(s^{r_1+r_2}).
            \end{align*}
        
        In our setting, $F$ is totally real of degree $n$ and $K$ is totally imaginary of degree $2n$. Since $K$ has no real embeddings and $F$ has no complex embeddings, $r_{1, F} + r_{2, F}= n = r_{1, K} + r_{2, K}.$ Therefore, both $\zeta_K(s)$ and $\zeta_F(s)$ have zeroes of order $n-1$ at $s=0$. Using $h_F = 1$, we have that
            \begin{align*}
                    \frac{\zeta_K(s)}{\zeta_F(s)} 
                    = \frac{2 h_K R_Ks^{n-1} + O(s^n)}{\omega_K R_Fs^{n-1} + O(s^n)}.
                \end{align*}
        Sending $s \rightarrow 0$, the higher order terms $O(s^n)$ tend to zero, and so we have that
            \begin{align*}
                h_K = \frac{w_K}{2} \cdot \frac{R_F}{R_K} \cdot L(0, \chi_{K/F}).
            \end{align*}
        Using the identity (see e.g. \citep[p. 406]{shintaniclnum})
            \begin{align*}
                \frac{1}{2} \frac{R_F}{R_K} = \frac{1}{[\oa_F^{\times}:\oa_F^{\times, +}][\oa_F^{\times, +}:N_{K/F}\oa_K^{\times}]},
            \end{align*}
        where $N_{K/F}\oa_K^{\times} \coloneqq \{N_{K/F}(x) | x \in \oa_K^{\times}\}$, and combining the work of Shintani \citep{shintaniclnum} and \citep{ddf, dasgupta} with Theorem \ref{thm: main thm}, we obtain
            \begin{align*}
                h_K &= \frac{1}{n} \cdot \frac{w_K}{ [\oa_F^{\times}:\oa_F^{\times, +}][\oa_F^{\times, +}:N_{K/F}\oa_K^{\times}]} \sum_{\substack{\tau \in S_{n-1} \\ w_\tau \neq 0}} w_\tau &\\
                &\hspace{+2.0cm} \times \left\{\sum_{m=1}^{p^n-1} (-1)^{m}\sum_{i=1}^{\#\left(\oa_F \cap R^\tau(p\oa_F)\right)} \sum_{\substack{\tiny{(l_1,\dots,l_n)\in\mathbb{Z}_{\geq0}^{n}}\\\sum_{j=1}^{n}l_{j}=n}} \prod_{k=1}^n \frac{B_{l_k}(\tilde{x}_{\tau, k}(i, m))}{l_k!} \tr_{F/\Q}\left(\prod_{k=1}^n f_{\tau,k}^{l_k-1}\right)\right\}.&
            \end{align*}

    \subsection{Proof of Corollary 1.3}
        Since $K$ is totally real of degree $2n$, we have that $r_{1, K} + r_{2, K} = 2n$ and hence the Taylor series expansion at $s = 0$ for the Dedekind zeta function associated to $K$ is given by
            \begin{align*}
                \zeta_K(s) = \frac{-h_KR_K}{w_K} s^{2n - 1} + O(s^{2n}).
            \end{align*}
        Since $F$ is also totally real and $[F:\Q]=n$, we have that $r_{1, F} + r_{2, F} = n$ and hence $\zeta_F(s)$ has a zero of order $n-1$ at $s=0$. It Taylor series expansion at $s=0$ is given by
            \begin{align*}
                \zeta_F(s) = \frac{-h_FR_F}{w_F} s^{n-1} + O(s^n).
            \end{align*}
        Since both $K$ and $F$ are totally real, $w_F = w_K = 2$. Recall that $h_F = 1$. Hence we have that
            \begin{align*}
                    L(s, \chi_{K/F}) = \frac{\zeta_K(s)}{\zeta_F(s)} = \frac{h_K R_K s^{2n-1} + O(s^{2n})}{R_Fs^{n-1} + O(s^n)}.
            \end{align*}
        Since $\zeta_K(s)$ has a zero of order $2n-1$ and $\zeta_F(s)$ has a zero of order $n-1$ at $s=0$, the $n^{th}$ derivative of $\zeta_K(s)/\zeta_F(s)$ evaluated at $s=0$ gives the ratio of the first nonzero terms in their Taylor series expansions. Therefore, the $n^{th}$ derivative of $L(s, \chi_{K/F})$ evaluated at $0$ is given by
            \begin{align*}
                h_K\cdot \frac{R_K}{R_F} = \frac{L^{(n)}(0, \chi_{K/F})}{n!}.
            \end{align*}
        By Theorem \ref{thm: main thm} combined with (\ref{charval--}), we have that
            \begin{align*}
                L(s, \chi_{K/F}) &= (-1)^n \cdot N(p\oa_F)^{-s} \sum_{\substack{\tau \in S_{n-1} \\ w_\tau \neq 0}} w_\tau \sum_{m=1}^{p^n-1} (-1)^m \sum_{i=1}^{\#\left(\oa_F \cap R^\tau(p\oa_F)\right)} \zeta(s, A^\tau, \tilde{\mathbf{x}}_{\tau}(i, m)),&
            \end{align*}
        and differentiating $n$ times yields
            \begin{align*}
                L^{(n)}(s, \chi_{K/F}) 
                & = N(p\oa_F)^{-s} \sum_{k=0}^n (-1)^k {n \choose k} \left(\ln(N(p\oa_F))\right)^{n-k} \\
                &\hspace{+3cm} \times \left\{ \sum_{\substack{\tau \in S_{n-1} \\ w_\tau \neq 0}} w_\tau \sum_{m=1}^{p^n-1} (-1)^m\sum_{i=1}^{\#\left(\oa_F \cap R^\tau(p\oa_F)\right)} \zeta^{(k)}(s, A^\tau, \tilde{\mathbf{x}}_{\tau}(i, m))\right\}.&
            \end{align*}
        Therefore,
            \begin{align*}
                h_K \cdot \frac{R_K}{R_F}
                & = \frac{1}{n!} \sum_{k=0}^n (-1)^k {n \choose k} \left(\ln(N(p\oa_F))\right)^{n-k} \sum_{\substack{\tau \in S_{n-1} \\ w_\tau \neq 0}} w_\tau \sum_{m=1}^{p^n-1} (-1)^m \sum_{i=1}^{\#\left(\oa_F \cap R^\tau(p\oa_F)\right)} \zeta^{(k)}\left(0, A^\tau, \tilde{\mathbf{x}}_{\tau}(i, m)\right).&
            \end{align*}

\section{Examples Illustrating Corollary \ref{cor: 1}} \label{section: examples}
    \begin{example}[1] Let $F = \Q(\zeta_7 + \zeta_7^{-1})$ and let $K = F(\sqrt{-3})$. Using \texttt{SageMath}, we compute that the totally positive unit group is generated by $\varepsilon_1 = (\zeta_7 + \zeta_7^{-1})^2$ and $\varepsilon_2 = (\zeta_7 + \zeta_7^{-1} + 1)^2$ and
            \begin{table}[H]
            \centering
                \begin{tabular}{|c|c|c|c|c|c|c|c|}
                    \hline
                    $\tau$ & $w_\tau$ & $f_{\tau, 1}$ &   $f_{\tau, 2}$ & $f_{\tau, 3}$ & $I_{\tau, 1}$ & $I_{\tau, 2}$ & $I_{\tau, 3}$ \\
                    \hline
                    $\mathrm{id}$ & 1 & 1 & $(\zeta_7 + \zeta_7^{-1})^2$ & $(\zeta_7 + \zeta_7^{-1})^2(\zeta_7 + \zeta_7^{-1}+1)^2$ & [0,1) & (0, 1] & [0,1)\\
                    \hline
                    $(12)$ & 1 & 1 & $(\zeta_7 + \zeta_7^{-1}+1)^2$ & $(\zeta_7 + \zeta_7^{-1})^2(\zeta_7 + \zeta_7^{-1}+1)^2$ & (0,1] & [0, 1) & (0,1] \\
                    \hline
                \end{tabular}
                
                \caption{The algebraic integers $f_{\tau, j}$ and associated intervals $I_{\tau, j}$ determining the boundary type of $R^\tau(3\mathcal{O}_F)$}
            \end{table}
        \noindent Since the associated weights are nonzero, both $\mathcal{B}_{F, \mathrm{id}} = \{1, f_{\mathrm{id}, 2}, f_{\mathrm{id}, 3}\}$ and $\mathcal{B}_{F, (12)}=\{1, f_{(12), 2}, f_{(12), 3}\}$ form a $\Q$-basis for $F$. We have that
            \begin{align*} 
            \centering
            R^\mathrm{id}(3\mathcal{O}_F) =&\left\{z \in \frac{1}{3}\mathcal{O}_F ~\bigg|~ z =\sum_{j=1}^3 t_{z, \mathrm{id}, j}f_{\mathrm{id}, j}, \quad \textbf{t}_{z, \mathrm{id}} \in [0,1) \times (0,1] \times [0,1) \right\}& \\ 
            R^{(12)}(3\mathcal{O}_F) =&\left\{z \in \frac{1}{3}\mathcal{O}_F ~\bigg|~ z =\sum_{j=1}^3 t_{z, (12), j}f_{(12), j}, \quad \textbf{t}_{z,(12)} \in (0, 1] \times [0, 1) \times (0,1] \right\}.&
            \end{align*}
        
        The irreducible polynomial for a primitive element $\theta_F = -(\zeta_7+\zeta_7^{-1})$ of $F$ is a primitive polynomial in $\mathbb{F}_{3}$, so we may take $\rho = \theta_F$ to be a generator of $\mathbb{F}_{27}^\times$. One checks that the minimal polynomial for $\rho$ is given by
            \begin{align*}
                h_\rho(x) = x^3 - x^2 - 2x + 1.
            \end{align*}
        We form the linear system
            \begin{align*}
                A_{F, \rho}(z)\textbf{X}_{F, \rho} = \begin{pmatrix}
                    1&0&z\\
                    -z&1&-2z\\
                    0&-z&1-z
                \end{pmatrix}\begin{pmatrix} X_{F, \rho, 1} \\ X_{F, \rho, 2} \\ X_{F, \rho, 3}\end{pmatrix} = \begin{pmatrix} -1 \\ 2 \\ 1\end{pmatrix},
            \end{align*}
        and by applying Cramer's Rule, we compute the unique rational functions solving this system to be
            \begin{align*}
                X_{F, \rho, 1} = \frac{-1}{1-z-2z^2+z^3}, \quad X_{F, \rho, 2} = \frac{2-z}{1-z-2z^2+z^3}, \quad \text{ and } \quad X_{F, \rho, 3} = \frac{1+2z-z^2}{1-z-2z^2+z^3}.
            \end{align*}
        For each $1 \leq m \leq 26$, we form the $3$-tuple $\textbf{x}(m) = (x_1(m), x_2(m), x_3(m))$, where $x_i(m)$ denotes the $m^{th}$ coefficient in the Taylor series expansion around $0$ of $X_{F, \rho, i}$ for $i = 1, 2, 3$, as can be seen in the second column of Table \ref{table:pts1}. The $3$-tuple $\overline{\textbf{x}}(m)$ is the minimal nonnegative representative of the equivalence class of $\textbf{x}(m)$ in $\mathbb{F}_3^3$ and gives the coordinates of $\rho^{3+m}$ in the basis $\mathcal{B}_{\mathbb{F}_{27}, \rho} = \{1, \rho, \rho^2\}$. These coordinate vectors are listed in column three of Table \ref{table:pts1}. Since $\rho = \theta_F$ is a primitive element of $\mathbb{F}_{27}$, the change of basis transformation $T_{\mathbb{F}_{27}, \mathcal{B}_{\rho}, \mathcal{B}_{\theta_F}} = I$. We apply the isomorphism $\varphi$ to $\rho^{3+m}$ given by sending $\theta_F$ to $\theta_F + 3\mathcal{O}_F$. Since we have that $\rho^{3+m} = \overline{x}_1(m) + \overline{x}_2(m)\theta_F + \overline{x}_3(m)\theta_F^2$ in $\mathbb{F}_{27}$, $\overline{\textbf{x}}(m)$ also gives the coordinates in the $\Q$-basis $\mathcal{B}_{\theta_F}$ for $F$ of a representative for the equivalence class of $\varphi(\rho^{3+m})$ in $\mathcal{O}_F/3\mathcal{O}_F$. We compute the change of basis transformations from $\mathcal{B}_{F, \theta_F}$ to $\mathcal{B}_{F, \tau}$ to be
            \begin{align*}
                T_\mathrm{id} = \begin{pmatrix}
                    1 & \frac{1}{3} & 0 \\
                    0 & \frac{2}{3} & 1 \\
                    0 & -\frac{1}{3} & 0
                \end{pmatrix} \quad \text{ and } \quad T_{(12)} = \begin{pmatrix}
                    1 & 1 & 1\\
                    0 & -2 & -3 \\
                    0 & 1 & 2
                \end{pmatrix}.
            \end{align*}
        Applying the transformation to $\overline{\textbf{x}}(m)$ and dividing by $3$, we obtain the vector of coordinates $\frac{1}{3}T_\tau\overline{\textbf{x}}(m)$ of the element $\frac{1}{3}\theta_F^{3+m}$ lying in the fractional ideal $\frac{1}{3}\mathcal{O}_F$ in the basis $\mathcal{B}_{F, \tau}$, as can be seen in the fourth column of Table \ref{table:pts1} and third column of Table \ref{table:pts2}. Recall that the Shintani set $R^\tau(3\mathcal{O}_F)$ forms a complete set of coset representatives for the group $G_\tau(3\mathcal{O}_F) = \frac{1}{3}\mathcal{O}_F/\bigoplus_{i=1}^3\Z f_{\tau, i}$. The unique corresponding point in the Shintani set $R^\tau(3\mathcal{O}_F)$ is computed by moding by $\bigoplus_{i=1}^3 \Z f_{\tau, i}$ to find the minimal representative for the equivalence class of $\frac{1}{3}\theta_F^{n+m}$ in $G_\tau(3\mathcal{O}_F)$ lying in $I_{\tau, 1} \times I_{\tau, 2} \times I_{\tau, 3}$. For each permutation $\tau$ in $S_2$, applying the modified fractional part functions associated to each interval $I_{\tau, i}$ to the $i^{th}$ coordinate of $\frac{1}{3}T_\tau\overline{\textbf{x}}(m)$ to obtain $\tilde{\textbf{x}}_\tau(1, m)$ corresponds to moding by $\bigoplus_{i=1}^3 \Z f_{\tau, i}$. Thus $\tilde{\textbf{x}}_\tau(1, m)$ gives the coordinates in the basis $\mathcal{B}_{F, \tau}$ of the unique coset representative $\vartheta(\rho^{3+m})$ of $\frac{1}{3}\theta_F^{3+m}$ in $R^\tau(3\mathcal{O}_F)$. Recall that the set of all such elements $C_\tau = \{\vartheta(\rho^{3+m}) : 1 \leq m \leq 26\},$ together with the identity element of $R^\tau(3\mathcal{O}_F)$, is a complete set of coset representatives for $R^\tau(3\mathcal{O}_F)/\mathrm{ker}(\pi_{3, \tau})$. We compute that
            \begin{align*}
                \mathrm{ker}(\pi_{3, \mathrm{id}}) = R^\mathrm{id}(3\mathcal{O}_F) \cap \mathcal{O}_F &= \left\{(0,1,0), \left(\frac{2}{3}, \frac{1}{3}, \frac{1}{3}\right), \left(\frac{1}{3}, \frac{2}{3}, \frac{2}{3}\right)\right\}&  \\
                \mathrm{ker}(\pi_{3, (12)}) = R^{(12)}(3\mathcal{O}_F) \cap \mathcal{O}_F &= \{(1,0,1)\}&.
            \end{align*}
        Therefore $R^\mathrm{id}(3\mathcal{O}_F)$ has $81$ elements, $R^{(12)}(3\mathcal{O}_F)$ has $27$ elements, and $\pi_{3, (12)}$ is an isomorphism between $R^{(12)}(3\mathcal{O}_F)$ and $\mathcal{O}_F/3\mathcal{O}_F$. We compute all nontrivial elements of $R^\mathrm{id}(3\mathcal{O}_F)$ by adding the vector of coordinates of the $i^{th}$ element $w(i)$ of $\mathrm{ker}(\pi_{3, \tau})$ in the basis $\mathcal{B}_{F, \tau}$ to each $\tilde{\textbf{x}}_\mathrm{id}(1, m)$ in $C_\mathrm{id}$ to obtain $\tilde{\textbf{x}}_\mathrm{id}(i, m)$ lying in the translated set by $C_\mathrm{id} \oplus w(i)$. Since the group law of $R^\mathrm{id}(3\mathcal{O}_F)$ is $\oplus$, the element $\vartheta(\rho^{3+m}) \oplus w(i)$ is computed by adding together the coordinates $\tilde{\textbf{x}}_\mathrm{id}(1, m)$ of 
        $\vartheta(\rho^{3+m})$ and those of $w(i)$ in the basis $\mathcal{B}_{F, \mathrm{id}}$ and computing the modified fractional part of each coordinate to obtain the unique coset representative in $R^\mathrm{id}(3\mathcal{O}_F)$. Thus the nontrivial elements of the Shintani sets are given by
            \begin{align*}
                R^\mathrm{id}(3\mathcal{O}_F) - \{(0,1,0)\} = C_{\mathrm{id}} \cup C_{\mathrm{id}} \oplus \left(\frac{2}{3}, \frac{1}{3}, \frac{1}{3}\right) \cup C_{\mathrm{id}}\oplus \left(\frac{1}{3}, \frac{2}{3}, \frac{2}{3}\right) \quad \text{ and } \quad R^{(12)}(3\mathcal{O}_F) -\{(1,0,1)\} = C_{(12)}.
            \end{align*}
        For $\tau = \mathrm{id}$, these sets are listed in the fifth, sixth, and seventh columns of Table \ref{table:pts1} and for $\tau = (12)$, it is listed in the fourth column of Table \ref{table:pts2}. Since the character value of an element lying in $\mathrm{ker}(\pi_{3, \tau})$ is $0$, we neglect these values in our calculations of the class number.
        
    \begin{table}[H]
    \centering
    \begin{tabular}
    {|c|c|c|c|c|c|c|}
    \hline
    &  &  & & \small$C_{\mathrm{id}}$ & \small$C_{\mathrm{id}} \oplus \big(\frac{2}{3}, \frac{1}{3},\frac{1}{3}\big)$&\small$C_{\mathrm{id}} \oplus \big(\frac{1}{3}, \frac{2}{3},\frac{2}{3}\big)$\\
    \textbf{m} & $\textbf{x}(m) = (x_1(m), x_2(m), x_2(m))$ & $\overline{\textbf{x}}(m)$ & $\frac{1}{3}T_\mathrm{id}\overline{\textbf{x}}(m)$ &  $\tilde{\textbf{x}}_{\mathrm{id}}(1,m)$&$\tilde{\textbf{x}}_{\mathrm{id}}(2,m)$&$\tilde{\textbf{x}}_{\mathrm{id}}(3,m)$\\
    \hline
    1 & (-1, 1, -3) & (2, 1, 0) & $\big(\frac{7}{9},\frac{2}{9},$-$\frac{1}{9}\big)$ & $\big(\frac{7}{9},\frac{2}{9},\frac{8}{9}\big)$ & $\big(\frac{4}{9},\frac{5}{9},\frac{2}{9}\big)$ & $\big(\frac{1}{9},\frac{8}{9},\frac{5}{9}\big)$\\ \hline 
    2 & (-3, 5, -4) & (0, 2, 1) & $\big(\frac{2}{9},\frac{7}{9},$-$\frac{2}{9}\big)$ & $\big(\frac{2}{9},\frac{7}{9},\frac{7}{9}\big)$ & $\big(\frac{8}{9},\frac{1}{9},\frac{1}{9}\big)$ & $\big(\frac{5}{9},\frac{4}{9},\frac{4}{9}\big)$\\ \hline 
    3 & (-4, 5, -9) & (2, 2, 0) & $\big(\frac{8}{9},\frac{4}{9},$-$\frac{2}{9}\big)$ & $\big(\frac{8}{9},\frac{4}{9},\frac{7}{9}\big)$ & $\big(\frac{5}{9},\frac{7}{9},\frac{1}{9}\big)$ & $\big(\frac{2}{9},\frac{1}{9},\frac{4}{9}\big)$\\ \hline 
    4 & (-9, 14, -14) & (0, 2, 2) & $\big(\frac{2}{9},\frac{10}{9},$-$\frac{2}{9}\big)$ & $\big(\frac{2}{9},\frac{1}{9},\frac{7}{9}\big)$ & $\big(\frac{8}{9},\frac{4}{9},\frac{1}{9}\big)$ & $\big(\frac{5}{9},\frac{7}{9},\frac{4}{9}\big)$\\ \hline 
    5 & (-14, 19, -28) & (1, 1, 1) & $\big(\frac{4}{9},\frac{5}{9},$-$\frac{1}{9}\big)$ & $\big(\frac{4}{9},\frac{5}{9},\frac{8}{9}\big)$ & $\big(\frac{1}{9},\frac{8}{9},\frac{2}{9}\big)$ & $\big(\frac{7}{9},\frac{2}{9},\frac{5}{9}\big)$\\ \hline 
    6 & (-28, 42, -47) & (2, 0, 2) & $\big(\frac{2}{3},\frac{2}{3},0\big)$ & $\big(\frac{2}{3},\frac{2}{3},0\big)$ & $\big(\frac{1}{3},1,\frac{1}{3}\big)$ & $\big(0,\frac{1}{3},\frac{2}{3}\big)$\\ \hline 
    7 & (-47, 66, -89) & (1, 0, 2) & $\big(\frac{1}{3},\frac{2}{3},0\big)$ & $\big(\frac{1}{3},\frac{2}{3},0\big)$ & $\big(0,1,\frac{1}{3}\big)$ & $\big(\frac{2}{3},\frac{1}{3},\frac{2}{3}\big)$\\ \hline 
    8 & (-89, 131, -155) & (1, 2, 2) & $\big(\frac{5}{9},\frac{10}{9},$-$\frac{2}{9}\big)$ & $\big(\frac{5}{9},\frac{1}{9},\frac{7}{9}\big)$ & $\big(\frac{2}{9},\frac{4}{9},\frac{1}{9}\big)$ & $\big(\frac{8}{9},\frac{7}{9},\frac{4}{9}\big)$\\ \hline 
    9 & (-155, 221, -286) & (1, 2, 1) & $\big(\frac{5}{9},\frac{7}{9},$-$\frac{2}{9}\big)$ & $\big(\frac{5}{9},\frac{7}{9},\frac{7}{9}\big)$ & $\big(\frac{2}{9},\frac{1}{9},\frac{1}{9}\big)$ & $\big(\frac{8}{9},\frac{4}{9},\frac{4}{9}\big)$\\ \hline 
    10 & (-286, 417, -507) & (2, 0, 0) & $\big(\frac{2}{3},0,0\big)$ & $\big(\frac{2}{3},1,0\big)$ & $\big(\frac{1}{3},\frac{1}{3},\frac{1}{3}\big)$ & $\big(0,\frac{2}{3},\frac{2}{3}\big)$\\ \hline 
    11 & (-507, 728, -924) & (0, 2, 0) & $\big(\frac{2}{9},\frac{4}{9},$-$\frac{2}{9}\big)$ & $\big(\frac{2}{9},\frac{4}{9},\frac{7}{9}\big)$ & $\big(\frac{8}{9},\frac{7}{9},\frac{1}{9}\big)$ & $\big(\frac{5}{9},\frac{1}{9},\frac{4}{9}\big)$\\ \hline 
    12 & (-924, 1341, -1652) & (0, 0, 2) & $\big(0,\frac{2}{3},0\big)$ & $\big(0,\frac{2}{3},0\big)$ & $\big(\frac{2}{3},1,\frac{1}{3}\big)$ & $\big(\frac{1}{3},\frac{1}{3},\frac{2}{3}\big)$\\ \hline 
    13 & (-1652, 2380, -2993) & (1, 1, 2) & $\big(\frac{4}{9},\frac{8}{9},$-$\frac{1}{9}\big)$ & $\big(\frac{4}{9},\frac{8}{9},\frac{8}{9}\big)$ & $\big(\frac{1}{9},\frac{2}{9},\frac{2}{9}\big)$ & $\big(\frac{7}{9},\frac{5}{9},\frac{5}{9}\big)$\\ \hline 
    14 & (-2993, 4334, -5373) & (1, 2, 0) & $\big(\frac{5}{9},\frac{4}{9},$-$\frac{2}{9}\big)$ & $\big(\frac{5}{9},\frac{4}{9},\frac{7}{9}\big)$ & $\big(\frac{2}{9},\frac{7}{9},\frac{1}{9}\big)$ & $\big(\frac{8}{9},\frac{1}{9},\frac{4}{9}\big)$\\ \hline 
    15 & (-5373, 7753, -9707) & (0, 1, 2) & $\big(\frac{1}{9},\frac{8}{9},$-$\frac{1}{9}\big)$ & $\big(\frac{1}{9},\frac{8}{9},\frac{8}{9}\big)$ & $\big(\frac{7}{9},\frac{2}{9},\frac{2}{9}\big)$ & $\big(\frac{4}{9},\frac{5}{9},\frac{5}{9}\big)$\\ \hline 
    16 & (-9707, 14041, -17460) & (1, 1, 0) & $\big(\frac{4}{9},\frac{2}{9},$-$\frac{1}{9}\big)$ & $\big(\frac{4}{9},\frac{2}{9},\frac{8}{9}\big)$ & $\big(\frac{1}{9},\frac{5}{9},\frac{2}{9}\big)$ & $\big(\frac{7}{9},\frac{8}{9},\frac{5}{9}\big)$\\ \hline 
    17 & (-17460, 25213, -31501) & (0, 1, 1) & $\big(\frac{1}{9},\frac{5}{9},$-$\frac{1}{9}\big)$ & $\big(\frac{1}{9},\frac{5}{9},\frac{8}{9}\big)$ & $\big(\frac{7}{9},\frac{8}{9},\frac{2}{9}\big)$ & $\big(\frac{4}{9},\frac{2}{9},\frac{5}{9}\big)$\\ \hline 
    18 & (-31501, 45542, -56714) & (2, 2, 2) & $\big(\frac{8}{9},\frac{10}{9},$-$\frac{2}{9}\big)$ & $\big(\frac{8}{9},\frac{1}{9},\frac{7}{9}\big)$ & $\big(\frac{5}{9},\frac{4}{9},\frac{1}{9}\big)$ & $\big(\frac{2}{9},\frac{7}{9},\frac{4}{9}\big)$\\ \hline 
    19 & (-56714, 81927, -102256) & (1, 0, 1) & $\big(\frac{1}{3},\frac{1}{3},0\big)$ & $\big(\frac{1}{3},\frac{1}{3},0\big)$ & $\big(0,\frac{2}{3},\frac{1}{3}\big)$ & $\big(\frac{2}{3},1,\frac{2}{3}\big)$\\ \hline 
    20 & (-102256, 147798, -184183) & (2, 0, 1) & $\big(\frac{2}{3},\frac{1}{3},0\big)$ & $\big(\frac{2}{3},\frac{1}{3},0\big)$ & $\big(\frac{1}{3},\frac{2}{3},\frac{1}{3}\big)$ & $\big(0,1,\frac{2}{3}\big)$\\ \hline 
    21 & (-184183, 266110, -331981) & (2, 1, 1) & $\big(\frac{7}{9},\frac{5}{9},$-$\frac{1}{9}\big)$ & $\big(\frac{7}{9},\frac{5}{9},\frac{8}{9}\big)$ & $\big(\frac{4}{9},\frac{8}{9},\frac{2}{9}\big)$ & $\big(\frac{1}{9},\frac{2}{9},\frac{5}{9}\big)$\\ \hline 
    22 & (-331981, 479779, -598091) & (2, 1, 2) & $\big(\frac{7}{9},\frac{8}{9},$-$\frac{1}{9}\big)$ & $\big(\frac{7}{9},\frac{8}{9},\frac{8}{9}\big)$ & $\big(\frac{4}{9},\frac{2}{9},\frac{2}{9}\big)$ & $\big(\frac{1}{9},\frac{5}{9},\frac{5}{9}\big)$\\ \hline 
    23 & (-598091, 864201, -1077870) & (1, 0, 0) & $\big(\frac{1}{3},0,0\big)$ & $\big(\frac{1}{3},1,0\big)$ & $\big(0,\frac{1}{3},\frac{1}{3}\big)$ & $\big(\frac{2}{3},\frac{2}{3},\frac{2}{3}\big)$\\ \hline 
    24 & (-1077870, 1557649, -1942071) & (0, 1, 0) & $\big(\frac{1}{9},\frac{2}{9},$-$\frac{1}{9}\big)$ & $\big(\frac{1}{9},\frac{2}{9},\frac{8}{9}\big)$ & $\big(\frac{7}{9},\frac{5}{9},\frac{2}{9}\big)$ & $\big(\frac{4}{9},\frac{8}{9},\frac{5}{9}\big)$\\ \hline 
    25 & (-1942071, 2806272, -3499720) & (0, 0, 1) & $\big(0,\frac{1}{3},0\big)$ & $\big(0,\frac{1}{3},0\big)$ & $\big(\frac{2}{3},\frac{2}{3},\frac{1}{3}\big)$ & $\big(\frac{1}{3},1,\frac{2}{3}\big)$\\ \hline 
    26 & (-3499720, 5057369, -6305992) & (2, 2, 1) & $\big(\frac{8}{9},\frac{7}{9},$-$\frac{2}{9}\big)$ & $\big(\frac{8}{9},\frac{7}{9},\frac{7}{9}\big)$ & $\big(\frac{5}{9},\frac{1}{9},\frac{1}{9}\big)$ & $\big(\frac{2}{9},\frac{4}{9},\frac{4}{9}\big)$\\ \hline
     \end{tabular}
            \caption{Computation of the nontrivial elements of $R^\mathrm{id}(3\mathcal{O}_F)$}
            \label{table:pts1}
    \end{table}

    \begin{table}[H]
    \centering
    \begin{tabular}
    { |c|c|c|c|c|}
    \hline
     &  & & \small$C_{(12)}$\\
    \textbf{m} & $\overline{\textbf{x}}(m)$ & $\frac{1}{3}T_{(12)}\overline{\textbf{x}}(m)$ &  $\tilde{\textbf{x}}_{(12)}(1, m)$\\
    \hline
    1 & (2, 1, 0) & $\big(1,$-$\frac{2}{3},\frac{1}{3}\big)$ & $\big(1,\frac{1}{3},\frac{1}{3}\big)$\\ \hline 
    2 & (0, 2, 1) & $\big(1,$-$\frac{7}{3},\frac{4}{3}\big)$ & $\big(1,\frac{2}{3},\frac{1}{3}\big)$\\ \hline 
    3 & (2, 2, 0) & $\big(\frac{4}{3},$-$\frac{4}{3},\frac{2}{3}\big)$ & $\big(\frac{1}{3},\frac{2}{3},\frac{2}{3}\big)$\\ \hline 
    4 & (0, 2, 2) & $\big(\frac{4}{3},$-$\frac{10}{3},2\big)$ & $\big(\frac{1}{3},\frac{2}{3},1\big)$\\ \hline 
    5 & (1, 1, 1) & $\big(1,$-$\frac{5}{3},1\big)$ & $\big(1,\frac{1}{3},1\big)$\\ \hline 
    6 & (2, 0, 2) & $\big(\frac{4}{3},$-$2,\frac{4}{3}\big)$ & $\big(\frac{1}{3},0,\frac{1}{3}\big)$\\ \hline 
    7 & (1, 0, 2) & $\big(1,$-$2,\frac{4}{3}\big)$ & $\big(1,0,\frac{1}{3}\big)$\\ \hline 
    8 & (1, 2, 2) & $\big(\frac{5}{3},$-$\frac{10}{3},2\big)$ & $\big(\frac{2}{3},\frac{2}{3},1\big)$\\ \hline 
    9 & (1, 2, 1) & $\big(\frac{4}{3},$-$\frac{7}{3},\frac{4}{3}\big)$ & $\big(\frac{1}{3},\frac{2}{3},\frac{1}{3}\big)$\\ \hline 
    10 & (2, 0, 0) & $\big(\frac{2}{3},0,0\big)$ & $\big(\frac{2}{3},0,1\big)$\\ \hline 
    11 & (0, 2, 0) & $\big(\frac{2}{3},$-$\frac{4}{3},\frac{2}{3}\big)$ & $\big(\frac{2}{3},\frac{2}{3},\frac{2}{3}\big)$\\ \hline 
    12 & (0, 0, 2) & $\big(\frac{2}{3},$-$2,\frac{4}{3}\big)$ & $\big(\frac{2}{3},0,\frac{1}{3}\big)$\\ \hline 
    13 & (1, 1, 2) & $\big(\frac{4}{3},$-$\frac{8}{3},\frac{5}{3}\big)$ & $\big(\frac{1}{3},\frac{1}{3},\frac{2}{3}\big)$\\ \hline 
    14 & (1, 2, 0) & $\big(1,$-$\frac{4}{3},\frac{2}{3}\big)$ & $\big(1,\frac{2}{3},\frac{2}{3}\big)$\\ \hline 
    15 & (0, 1, 2) & $\big(1,$-$\frac{8}{3},\frac{5}{3}\big)$ & $\big(1,\frac{1}{3},\frac{2}{3}\big)$\\ \hline 
    16 & (1, 1, 0) & $\big(\frac{2}{3},$-$\frac{2}{3},\frac{1}{3}\big)$ & $\big(\frac{2}{3},\frac{1}{3},\frac{1}{3}\big)$\\ \hline 
    17 & (0, 1, 1) & $\big(\frac{2}{3},$-$\frac{5}{3},1\big)$ & $\big(\frac{2}{3},\frac{1}{3},1\big)$\\ \hline 
    18 & (2, 2, 2) & $\big(2,$-$\frac{10}{3},2\big)$ & $\big(1,\frac{2}{3},1\big)$\\ \hline 
    19 & (1, 0, 1) & $\big(\frac{2}{3},$-$1,\frac{2}{3}\big)$ & $\big(\frac{2}{3},0,\frac{2}{3}\big)$\\ \hline 
    20 & (2, 0, 1) & $\big(1,$-$1,\frac{2}{3}\big)$ & $\big(1,0,\frac{2}{3}\big)$\\ \hline 
    21 & (2, 1, 1) & $\big(\frac{4}{3},$-$\frac{5}{3},1\big)$ & $\big(\frac{1}{3},\frac{1}{3},1\big)$\\ \hline 
    22 & (2, 1, 2) & $\big(\frac{5}{3},$-$\frac{8}{3},\frac{5}{3}\big)$ & $\big(\frac{2}{3},\frac{1}{3},\frac{2}{3}\big)$\\ \hline 
    23 & (1, 0, 0) & $\big(\frac{1}{3},0,0\big)$ & $\big(\frac{1}{3},0,1\big)$\\ \hline 
    24 & (0, 1, 0) & $\big(\frac{1}{3},$-$\frac{2}{3},\frac{1}{3}\big)$ & $\big(\frac{1}{3},\frac{1}{3},\frac{1}{3}\big)$\\ \hline 
    25 & (0, 0, 1) & $\big(\frac{1}{3},$-$1,\frac{2}{3}\big)$ & $\big(\frac{1}{3},0,\frac{2}{3}\big)$\\ \hline 
    26 & (2, 2, 1) & $\big(\frac{5}{3},$-$\frac{7}{3},\frac{4}{3}\big)$ & $\big(\frac{2}{3},\frac{2}{3},\frac{1}{3}\big)$\\ \hline    
     \end{tabular}
            \caption{Computation of the nontrivial elements of $R^{(12)}(3\mathcal{O}_F)$}
            \label{table:pts2}
    \end{table}
    
    \begin{table}[H]
        \centering
        \begin{tabular}{|c|c|c|c||c|c|}
            \hline
            \multirow{2}{*}{\textbf{m}} & \multicolumn{3}{c||}{$S_\mathrm{id}(i,m)$} & $S_{(12)}(i,m)$ & $\mathbf{\Sigma}$ \\
            & $1$ & $2$ & $3$ & $1$ & \\
            \hline
            1 & $\frac{109}{108}$ & -$\frac{65}{108}$ & -$\frac{53}{108}$ & $\frac{7}{36}$ & $\frac{1}{9}$ \\ \hline 
            2 & $\frac{65}{324}$ & -$\frac{133}{324}$ & -$\frac{7}{324}$ & $\frac{13}{324}$ & -$\frac{31}{162}$ \\ \hline 
            3 & $\frac{317}{324}$ & -$\frac{7}{324}$ & -$\frac{313}{324}$ & -$\frac{5}{108}$ & -$\frac{1}{18}$ \\ \hline 
            4 & $\frac{161}{324}$ & -$\frac{91}{324}$ & -$\frac{145}{324}$ & $\frac{41}{324}$ & -$\frac{17}{162}$ \\ \hline 
            5 & $\frac{91}{324}$ & -$\frac{323}{324}$ & $\frac{217}{324}$ & $\frac{71}{324}$ & $\frac{14}{81}$ \\ \hline 
            6 & -$\frac{103}{324}$ & $\frac{113}{324}$ & -$\frac{49}{324}$ & $\frac{215}{324}$ & $\frac{44}{81}$ \\ \hline 
            7 & -$\frac{217}{324}$ & -$\frac{217}{324}$ & $\frac{233}{324}$ & $\frac{97}{324}$ & -$\frac{26}{81}$ \\ \hline 
            8 & -$\frac{25}{108}$ & $\frac{107}{108}$ & -$\frac{61}{108}$ & -$\frac{121}{324}$ & -$\frac{29}{162}$ \\ \hline 
            9 & $\frac{89}{108}$ & -$\frac{61}{108}$ & $\frac{65}{108}$ & -$\frac{133}{324}$ & $\frac{73}{162}$ \\ \hline 
            10 & -$\frac{139}{324}$ & $\frac{329}{324}$ & $\frac{23}{324}$ & -$\frac{11}{36}$ & $\frac{19}{54}$ \\ \hline 
            11 & -$\frac{119}{324}$ & $\frac{61}{324}$ & -$\frac{119}{324}$ & $\frac{25}{36}$ & $\frac{4}{27}$ \\ \hline 
            12 & $\frac{11}{108}$ & -$\frac{73}{108}$ & $\frac{47}{108}$ & $\frac{65}{324}$ & $\frac{5}{81}$ \\ \hline 
            13 & $\frac{103}{324}$ & -$\frac{203}{324}$ & $\frac{301}{324}$ & -$\frac{125}{324}$ & $\frac{19}{81}$ \\ \hline 
            14 & -$\frac{65}{108}$ & $\frac{109}{108}$ & -$\frac{53}{108}$ & -$\frac{1}{36}$ & -$\frac{1}{9}$ \\ \hline 
            15 & -$\frac{133}{324}$ & $\frac{65}{324}$ & -$\frac{7}{324}$ & $\frac{139}{324}$ & $\frac{16}{81}$ \\ \hline 
            16 & -$\frac{7}{324}$ & $\frac{317}{324}$ & -$\frac{313}{324}$ & -$\frac{5}{108}$ & -$\frac{1}{18}$ \\ \hline 
            17 & -$\frac{91}{324}$ & $\frac{161}{324}$ & -$\frac{145}{324}$ & $\frac{167}{324}$ & $\frac{23}{81}$ \\ \hline 
            18 & -$\frac{323}{324}$ & $\frac{91}{324}$ & $\frac{217}{324}$ & $\frac{233}{324}$ & $\frac{109}{162}$ \\ \hline 
            19 & -$\frac{175}{324}$ & -$\frac{175}{324}$ & $\frac{329}{324}$ & -$\frac{1}{324}$ & -$\frac{11}{162}$ \\ \hline 
            20 & -$\frac{91}{324}$ & $\frac{233}{324}$ & $\frac{71}{324}$ & -$\frac{155}{324}$ & $\frac{29}{162}$ \\ \hline 
            21 & $\frac{107}{108}$ & -$\frac{25}{108}$ & -$\frac{61}{108}$ & -$\frac{49}{324}$ & $\frac{7}{162}$ \\ \hline 
            22 & -$\frac{61}{108}$ & $\frac{89}{108}$ & $\frac{65}{108}$ & -$\frac{133}{324}$ & $\frac{73}{162}$ \\ \hline 
            23 & -$\frac{175}{324}$ & -$\frac{49}{324}$ & $\frac{329}{324}$ & -$\frac{11}{36}$ & $\frac{1}{54}$ \\ \hline 
            24 & $\frac{61}{324}$ & -$\frac{119}{324}$ & -$\frac{119}{324}$ & $\frac{25}{36}$ & $\frac{4}{27}$ \\ \hline 
            25 & $\frac{89}{108}$ & $\frac{47}{108}$ & $\frac{17}{108}$ & -$\frac{205}{324}$ & $\frac{127}{162}$ \\ \hline 
            26 & -$\frac{203}{324}$ & $\frac{103}{324}$ & $\frac{301}{324}$ & -$\frac{125}{324}$ & $\frac{19}{81}$ \\ \hline
            $\mathbf{\Sigma}$ & -$\frac{2}{3}$ & $\frac{5}{3}$ & $\frac{13}{6}$ & $\frac{5}{6}$ & $\mathbf{4}$ \\
            \hline
        \end{tabular}
        \caption{Computation of $S_\tau(i,m) = (-1)^m \sum_{\substack{\tiny{(l_1,1_2,l_3)\in\mathbb{Z}_{\geq0}^{3}}\\1_1 + l_2+l_3=3}} \prod_{k=1}^3 \frac{B_{l_k}(\tilde{x}_{\tau,k}(i, m))}{l_k!} \,\tr_{F/\Q}\bigg(\prod_{k=1}^3 f_{\tau,k}^{l_k-1}\bigg)$}
        \label{table:sum}
    \end{table}
    \newpage
    \noindent Using \texttt{SageMath}, we compute $w_K=6$, $[\oa_F^{\times}:\oa_F^{\times, +}]=8$, and $[\oa_F^{\times, +}:N_{K/F}\oa_K^{\times}] =1.$
    Hence, adding the contributions from each $S_\tau(i,m)$ listed in Table \ref{table:sum}, we compute that the class number is given by
        \begin{align*}
        h_K=\frac{1}{3}\cdot\frac{6}{8}\left(-\frac{3}{2}+\frac{5}{3}+\frac{13}{6}+\frac{5}{6}\right) = 1.
        \end{align*}
    \end{example}
    \begin{example}[2] Let $F = \Q[x]/(x^3 - x^2 - 6x + 7) = \Q(\theta_F)$ where $\theta_F$ is a root of $x^3-x^2-6x+7$ and let $K = F(\sqrt{-3})$. Using \texttt{SageMath}, we compute that the totally positive unit group is generated by $\varepsilon_1 = 2\theta_F^2 + 3\theta_F -5$ and $\varepsilon_2 = \theta_F^2-2\theta_F+1$ and
            \begin{table}[H]
            \centering
                \begin{tabular}{|c|c|c|c|c|c|c|c|}
                    \hline
                    $\tau$ & $w_\tau$ & $f_{\tau, 1}$ &   $f_{\tau, 2}$ & $f_{\tau, 3}$ & $I_{\tau, 1}$ & $I_{\tau, 2}$ & $I_{\tau, 3}$ \\
                    \hline
                    $\mathrm{id}$ & 1 & 1 & $2\theta_F^2 + 3\theta_F -5$ & $(2\theta_F^2 + 3\theta_F -5)(\theta_F^2-2\theta_F+1)$ & [0,1) & (0, 1] & (0,1]\\
                    \hline
                    $(12)$ & 1 & 1 & $\theta_F^2-2\theta_F+1$ & $(2\theta_F^2 + 3\theta_F -5)(\theta_F^2-2\theta_F+1)$ & [0,1) & [0, 1) & (0,1] \\
                    \hline
                \end{tabular}
                
                \caption{The algebraic integers $f_{\tau, j}$ and associated intervals $I_{\tau, j}$ determining the boundary type of $R^\tau(3\mathcal{O}_F)$}
            \end{table}
        \noindent The irreducible polynomial for $\theta_F$ is a primitive polynomial in $\mathbb{F}_{3}$, so we may take $\rho = \theta_F$ to be a generator of $\mathbb{F}_{27}^\times$ and therefore
            \begin{align*}
                h_\rho(x) = x^3 - x^2 - 6x + 7.
            \end{align*}
        We form the linear system
            \begin{align*}
                A_{F, \rho}(z)\textbf{X}_{F, \rho} = \begin{pmatrix}
                    1&0&7z\\
                    -z&1&-6z\\
                    0&-z&1-z
                \end{pmatrix}\begin{pmatrix} X_{F, \rho, 1} \\ X_{F, \rho, 2} \\ X_{F, \rho, 3}\end{pmatrix} = \begin{pmatrix} -7 \\ 6 \\ 1\end{pmatrix},
            \end{align*}
        and by applying Cramer's Rule, we compute the unique rational functions solving this system to be
            \begin{align*}
                X_{F, \rho, 1} = \frac{-7}{1-6z-z^2+7z^3}, \quad X_{F, \rho, 2} = \frac{6-7z}{1-6z-z^2+7z^3}, \quad \text{ and } \quad X_{F, \rho, 3} = \frac{1+6z-7z^2}{1-6z-z^2+7z^3}.
            \end{align*}
        Since $\rho = \theta_F$ is a primitive element in $\mathbb{F}_{27}$, the change of basis transformation $T_{\mathbb{F}_{27}, \mathcal{B}_{\rho}, \mathcal{B}_{\theta_F}} = I$. We compute the change of basis transformations from $\mathcal{B}_{F, \theta_F}$ to $\mathcal{B}_{F, \tau}$ to be
            \begin{align*}
                T_\mathrm{id} = \begin{pmatrix}
                    1 & -2 & \frac{11}{2} \\
                    0 & 2 & -\frac{5}{2} \\
                    0 & -1 & \frac{3}{2}
                \end{pmatrix} \quad \text{ and } \quad T_{(12)} = \begin{pmatrix}
                    1 & \frac{16}{13} & \frac{19}{13}\\
                    0 & -\frac{4}{13} & \frac{5}{13} \\
                    0 & \frac{1}{13} & \frac{2}{13}
                \end{pmatrix}.
            \end{align*}
       We compute that
            \begin{align*}
                \mathrm{ker}(\pi_{3, \mathrm{id}}) = \left\{(1,0,0), \left(\frac{1}{2}, \frac{1}{2}, \frac{1}{2}\right)\right\}
            \end{align*}
            \begin{align*}
                &\mathrm{ker}(\pi_{3, (12)}) = \bigg\{(1,1,0), 
                \left(\frac{3}{13}, \frac{9}{13}, \frac{1}{13}\right),
                \left(\frac{12}{13}, \frac{10}{13}, \frac{4}{13}\right),
                \left(\frac{8}{13}, \frac{11}{13}, \frac{7}{13}\right),
                \left(\frac{4}{13}, \frac{12}{13}, \frac{10}{13}\right),\left(\frac{6}{13}, \frac{5}{13}, \frac{2}{13}\right),
                \left(\frac{2}{13}, \frac{6}{13}, \frac{5}{13}\right),& \\
                &\hspace{2cm}
               \left(\frac{11}{13}, \frac{7}{13}, \frac{8}{13}\right),
                \left(\frac{7}{13}, \frac{8}{13}, \frac{11}{13}\right),
                \left(\frac{9}{13}, \frac{1}{13}, \frac{3}{13}\right),
                \left(\frac{5}{13}, \frac{2}{13}, \frac{6}{13}\right),
                \left(\frac{1}{13}, \frac{3}{13}, \frac{9}{13}\right),
                \left(\frac{10}{13}, \frac{4}{13}, \frac{12}{13}\right)\bigg\}.&
            \end{align*}
        Therefore $R^\mathrm{id}(3\mathcal{O}_F)$ has $54$ elements, $R^{(12)}(3\mathcal{O}_F)$ has $351$ elements. We list the sets $C_\mathrm{id}$ and $C_{(12)}$ in columns five and seven of Table (\ref{ptsex2}). 
        \begin{table}[H]
            \centering
            \begin{tabular}
                {|c|c|c|c|c||c|c|}
                \hline
                &  &  &  &  $C_{\mathrm{id}}$& & $C_{(12)}$\\
                \textbf{m} & $\textbf{x}(m) = (x_1(m), x_2(m), x_2(m))$ & $\overline{\textbf{x}}(m)$ & $\frac{1}{3}T_\mathrm{id}\overline{\textbf{x}}(m)$ &  $\tilde{\textbf{x}}_{\mathrm{id}}(1,m)$&$\frac{1}{3}T_{(12)}\overline{\textbf{x}}(m)$&$\tilde{\textbf{x}}_{(12)}(1,m)$\\
                \hline
                1 & (-7, -1, 7) & (2, 2, 1) & $\big(\frac{7}{6},\frac{1}{2},$-$\frac{1}{6}\big)$ & $\big(\frac{1}{6},\frac{1}{2},\frac{5}{6}\big)$ & $\big(\frac{77}{39},$-$\frac{1}{13},\frac{4}{39}\big)$ & $\big(\frac{38}{39},\frac{12}{13},\frac{4}{39}\big)$\\ \hline 
                2 & (-49, 35, 6) & (2, 2, 0) & $\big($-$\frac{2}{3},\frac{4}{3},$-$\frac{2}{3}\big)$ & $\big(\frac{1}{3},\frac{1}{3},\frac{1}{3}\big)$ & $\big(\frac{58}{39},$-$\frac{8}{39},\frac{2}{39}\big)$ & $\big(\frac{19}{39},\frac{31}{39},\frac{2}{39}\big)$\\ \hline 
                3 & (-42, -13, 41) & (0, 2, 2) & $\big(\frac{7}{3},$-$\frac{1}{3},\frac{1}{3}\big)$ & $\big(\frac{1}{3},\frac{2}{3},\frac{1}{3}\big)$ & $\big(\frac{70}{39},\frac{2}{39},\frac{2}{13}\big)$ & $\big(\frac{31}{39},\frac{2}{39},\frac{2}{13}\big)$\\ \hline 
                4 & (-287, 204, 28) & (1, 0, 1) & $\big(\frac{13}{6},$-$\frac{5}{6},\frac{1}{2}\big)$ & $\big(\frac{1}{6},\frac{1}{6},\frac{1}{2}\big)$ & $\big(\frac{32}{39},\frac{5}{39},\frac{2}{39}\big)$ & $\big(\frac{32}{39},\frac{5}{39},\frac{2}{39}\big)$\\ \hline 
                5 & (-196, -119, 232) & (2, 1, 1) & $\big(\frac{11}{6},$-$\frac{1}{6},\frac{1}{6}\big)$ & $\big(\frac{5}{6},\frac{5}{6},\frac{1}{6}\big)$ & $\big(\frac{61}{39},\frac{1}{39},\frac{1}{13}\big)$ & $\big(\frac{22}{39},\frac{1}{39},\frac{1}{13}\big)$\\ \hline 
                6 & (-1624, 1196, 113) & (2, 2, 2) & $\big(3,$-$\frac{1}{3},\frac{1}{3}\big)$ & $\big(1,\frac{2}{3},\frac{1}{3}\big)$ & $\big(\frac{32}{13},\frac{2}{39},\frac{2}{13}\big)$ & $\big(\frac{6}{13},\frac{2}{39},\frac{2}{13}\big)$\\ \hline 
                7 & (-791, -946, 1309) & (1, 2, 1) & $\big(\frac{5}{6},\frac{1}{2},$-$\frac{1}{6}\big)$ & $\big(\frac{5}{6},\frac{1}{2},\frac{5}{6}\big)$ & $\big(\frac{64}{39},$-$\frac{1}{13},\frac{4}{39}\big)$ & $\big(\frac{25}{39},\frac{12}{13},\frac{4}{39}\big)$\\ \hline 
                8 & (-9163, 7063, 363) & (2, 1, 0) & $\big(0,\frac{2}{3},$-$\frac{1}{3}\big)$ & $\big(1,\frac{2}{3},\frac{2}{3}\big)$ & $\big(\frac{14}{13},$-$\frac{4}{39},\frac{1}{39}\big)$ & $\big(\frac{1}{13},\frac{35}{39},\frac{1}{39}\big)$\\ \hline 
                9 & (-2541, -6985, 7426) & (0, 2, 1) & $\big(\frac{1}{2},\frac{1}{2},$-$\frac{1}{6}\big)$ & $\big(\frac{1}{2},\frac{1}{2},\frac{5}{6}\big)$ & $\big(\frac{17}{13},$-$\frac{1}{13},\frac{4}{39}\big)$ & $\big(\frac{4}{13},\frac{12}{13},\frac{4}{39}\big)$\\ \hline 
                10 & (-51982, 42015, 441) & (2, 0, 0) & $\big(\frac{2}{3},0,0\big)$ & $\big(\frac{2}{3},0,0\big)$ & $\big(\frac{2}{3},0,0\big)$ & $\big(\frac{2}{3},1,0\big)$\\ \hline 
                11 & (-3087, -49336, 42456) & (0, 2, 0) & $\big($-$\frac{4}{3},\frac{4}{3},$-$\frac{2}{3}\big)$ & $\big(\frac{2}{3},\frac{1}{3},\frac{1}{3}\big)$ & $\big(\frac{32}{39},$-$\frac{8}{39},\frac{2}{39}\big)$ & $\big(\frac{32}{39},\frac{31}{39},\frac{2}{39}\big)$\\ \hline 
                12 & (-297192, 251649, -6880) & (0, 0, 2) & $\big(\frac{11}{3},$-$\frac{5}{3},1\big)$ & $\big(\frac{2}{3},\frac{1}{3},0\big)$ & $\big(\frac{38}{39},\frac{10}{39},\frac{4}{39}\big)$ & $\big(\frac{38}{39},\frac{10}{39},\frac{4}{39}\big)$\\ \hline 
                13 & (48160, -338472, 244769) & (1, 0, 2) & $\big(4,$-$\frac{5}{3},1\big)$ & $\big(1,\frac{1}{3},0\big)$ & $\big(\frac{17}{13},\frac{10}{39},\frac{4}{39}\big)$ & $\big(\frac{4}{13},\frac{10}{39},\frac{4}{39}\big)$\\ \hline 
                14 & (-1713383, 1516774, -93703) & (1, 1, 2) & $\big(\frac{10}{3},$-$1,\frac{2}{3}\big)$ & $\big(\frac{1}{3},0,\frac{2}{3}\big)$ & $\big(\frac{67}{39},\frac{2}{13},\frac{5}{39}\big)$ & $\big(\frac{28}{39},\frac{2}{13},\frac{5}{39}\big)$\\ \hline 
                15 & (655921, -2275601, 1423071) & (1, 1, 0) & $\big($-$\frac{1}{3},\frac{2}{3},$-$\frac{1}{3}\big)$ & $\big(\frac{2}{3},\frac{2}{3},\frac{2}{3}\big)$ & $\big(\frac{29}{39},$-$\frac{4}{39},\frac{1}{39}\big)$ & $\big(\frac{29}{39},\frac{35}{39},\frac{1}{39}\big)$\\ \hline 
                16 & (-9961497, 9194347, -852530) & (0, 1, 1) & $\big(\frac{7}{6},$-$\frac{1}{6},\frac{1}{6}\big)$ & $\big(\frac{1}{6},\frac{5}{6},\frac{1}{6}\big)$ & $\big(\frac{35}{39},\frac{1}{39},\frac{1}{13}\big)$ & $\big(\frac{35}{39},\frac{1}{39},\frac{1}{13}\big)$\\ \hline 
                17 & (5967710, -15076677, 8341817) & (2, 0, 2) & $\big(\frac{13}{3},$-$\frac{5}{3},1\big)$ & $\big(\frac{1}{3},\frac{1}{3},0\big)$ & $\big(\frac{64}{39},\frac{10}{39},\frac{4}{39}\big)$ & $\big(\frac{25}{39},\frac{10}{39},\frac{4}{39}\big)$\\ \hline 
                18 & (-58392719, 56018612, -6734860) & (1, 2, 2) & $\big(\frac{8}{3},$-$\frac{1}{3},\frac{1}{3}\big)$ & $\big(\frac{2}{3},\frac{2}{3},\frac{1}{3}\big)$ & $\big(\frac{83}{39},\frac{2}{39},\frac{2}{13}\big)$ & $\big(\frac{5}{39},\frac{2}{39},\frac{2}{13}\big)$\\ \hline 
                19 & (47144020, -98801879, 49283752) & (1, 1, 1) & $\big(\frac{3}{2},$-$\frac{1}{6},\frac{1}{6}\big)$ & $\big(\frac{1}{2},\frac{5}{6},\frac{1}{6}\big)$ & $\big(\frac{16}{13},\frac{1}{39},\frac{1}{13}\big)$ & $\big(\frac{3}{13},\frac{1}{39},\frac{1}{13}\big)$\\ \hline 
                20 & (-344986264, 342846532, -49518127) & (2, 1, 2) & $\big(\frac{11}{3},$-$1,\frac{2}{3}\big)$ & $\big(\frac{2}{3},0,\frac{2}{3}\big)$ & $\big(\frac{80}{39},\frac{2}{13},\frac{5}{39}\big)$ & $\big(\frac{2}{39},\frac{2}{13},\frac{5}{39}\big)$\\ \hline 
                21 & (346626889, -642095026, 293328405) & (1, 2, 0) & $\big($-$1,\frac{4}{3},$-$\frac{2}{3}\big)$ & $\big(1,\frac{1}{3},\frac{1}{3}\big)$ & $\big(\frac{15}{13},$-$\frac{8}{39},\frac{2}{39}\big)$ & $\big(\frac{2}{13},\frac{31}{39},\frac{2}{39}\big)$\\ \hline 
                22 & (-2053298835, 2106597319, -348766621) & (0, 1, 2) & $\big(3,$-$1,\frac{2}{3}\big)$ & $\big(1,0,\frac{2}{3}\big)$ & $\big(\frac{18}{13},\frac{2}{13},\frac{5}{39}\big)$ & $\big(\frac{5}{13},\frac{2}{13},\frac{5}{39}\big)$\\ \hline 
                23 & (2441366347, -4145898561, 1757830698) & (1, 0, 0) & $\big(\frac{1}{3},0,0\big)$ & $\big(\frac{1}{3},0,0\big)$ & $\big(\frac{1}{3},0,0\big)$ & $\big(\frac{1}{3},1,0\big)$\\ \hline 
                24 & (-12304814886, 12988350535, -2388067863) & (0, 1, 0) & $\big($-$\frac{2}{3},\frac{2}{3},$-$\frac{1}{3}\big)$ & $\big(\frac{1}{3},\frac{2}{3},\frac{2}{3}\big)$ & $\big(\frac{16}{39},$-$\frac{4}{39},\frac{1}{39}\big)$ & $\big(\frac{16}{39},\frac{35}{39},\frac{1}{39}\big)$\\ \hline 
                25 & (16716475041, -26633222064, 10600282672) & (0, 0, 1) & $\big(\frac{11}{6},$-$\frac{5}{6},\frac{1}{2}\big)$ & $\big(\frac{5}{6},\frac{1}{6},\frac{1}{2}\big)$ & $\big(\frac{19}{39},\frac{5}{39},\frac{2}{39}\big)$ & $\big(\frac{19}{39},\frac{5}{39},\frac{2}{39}\big)$\\ \hline 
                26 & (-74201978704, 80318171073, -16032939392) & (2, 0, 1) & $\big(\frac{5}{2},$-$\frac{5}{6},\frac{1}{2}\big)$ & $\big(\frac{1}{2},\frac{1}{6},\frac{1}{2}\big)$ & $\big(\frac{15}{13},\frac{5}{39},\frac{2}{39}\big)$ & $\big(\frac{2}{13},\frac{5}{39},\frac{2}{39}\big)$\\ \hline
            \end{tabular}
            \caption{Computation of the sets $C_{\mathrm{id}}$ and $C_{(12)}$ using the rational generating functions}
            \label{ptsex2}
        \end{table}
        
        To obtain $\tilde{\textbf{x}}_\mathrm{id}(2, m)$, $\left(\frac{1}{2},\frac{1}{2},\frac{1}{2}\right)$ is added to the vector $\tilde{\textbf{x}}_\mathrm{id}(1,m)$, and the result is reduced to lie in $R^\mathrm{id}(3\mathcal{O}_F)$ using the fractional part functions for the intervals $I_{\mathrm{id}, j}$, $j=1,2,3$. Similarly, to obtain $\tilde{\textbf{x}}_{(12)}(i, m)$ for each of the 12 nontrivial elements of $\mathrm{ker}(\pi_{3,(12)})$, we compute the sum of its coordinate vector and $\tilde{\textbf{x}}_{(12)}(1, m)$, and reduce the result to lie in $R^{(12)}(3\mathcal{O}_F)$ using the fractional part functions for the intervals $I_{(12), j}$, $j=1,2,3$.
        \begin{table}[H]
            \centering
            \begin{tabular}
                {|c|c|c|}
                \hline
                \multirow{2}{*}{\textbf{m}} & \multicolumn{2}{c|}{$S_\mathrm{id}(i,m)$} \\
                & $1$ & $2$ \\
                \hline
                1 & -$\frac{7}{54}$ & -$\frac{71}{108}$\\ \hline 
                2 & $\frac{133}{108}$ & -$\frac{19}{54}$\\ \hline 
                3 & -$\frac{139}{324}$ & $\frac{43}{162}$\\ \hline 
                4 & $\frac{155}{162}$ & $\frac{85}{324}$\\ \hline 
                5 & $\frac{43}{162}$ & -$\frac{139}{324}$\\ \hline 
                6 & -$\frac{71}{108}$ & -$\frac{7}{54}$\\ \hline 
                7 & $\frac{155}{162}$ & -$\frac{203}{324}$\\ \hline 
                8 & -$\frac{203}{324}$ & $\frac{155}{162}$\\ \hline 
                9 & $\frac{65}{81}$ & -$\frac{91}{324}$\\ \hline 
                10 & -$\frac{199}{324}$ & $\frac{65}{81}$\\ \hline 
                11 & -$\frac{139}{324}$ & $\frac{43}{162}$\\ \hline 
                12 & $\frac{43}{108}$ & -$\frac{7}{54}$\\ \hline 
                13 & $\frac{161}{324}$ & $\frac{65}{81}$\\ \hline 
                14 & $\frac{43}{108}$ & -$\frac{7}{54}$\\ \hline 
                15 & $\frac{133}{108}$ & -$\frac{19}{54}$\\ \hline 
                16 & $\frac{43}{162}$ & -$\frac{139}{324}$\\ \hline 
                17 & -$\frac{203}{324}$ & $\frac{155}{162}$\\ \hline 
                18 & -$\frac{139}{324}$ & $\frac{43}{162}$\\ \hline 
                19 & -$\frac{7}{54}$ & $\frac{43}{108}$\\ \hline 
                20 & $\frac{85}{324}$ & $\frac{155}{162}$\\ \hline 
                21 & $\frac{85}{324}$ & $\frac{155}{162}$\\ \hline 
                22 & $\frac{161}{324}$ & $\frac{65}{81}$\\ \hline 
                23 & $\frac{431}{324}$ & $\frac{65}{81}$\\ \hline 
                24 & -$\frac{139}{324}$ & $\frac{43}{162}$\\ \hline 
                25 & -$\frac{7}{54}$ & -$\frac{71}{108}$\\ \hline 
                26 & $\frac{65}{81}$ & -$\frac{91}{324}$\\ \hline 
                $\mathbf{\Sigma}_\mathrm{id}$ & $\frac{199}{36}$ & $\frac{155}{36}$ \\
                \hline
            \end{tabular}
            \caption{Computation of $S_\mathrm{id}(i,m) = (-1)^m \sum_{\substack{\tiny{(l_1,1_2,l_3)\in\mathbb{Z}_{\geq0}^{3}}\\1_1 + l_2+l_3=3}} \prod_{k=1}^3 \frac{B_{l_k}(\tilde{x}_{\mathrm{id},k}(i, m))}{l_k!} \,\tr_{F/\Q}\bigg(\prod_{k=1}^3 f_{\mathrm{id},k}^{l_k-1}\bigg)$}
            \label{table:sums1}
        \end{table}

        \begin{table}[H]
            \centering
            \begin{tabular}
                {|c|c|c|c|c|c|c|c|c|c|c|c|c|c|}
                \hline
                \multirow{2}{*}{\textbf{m}} & \multicolumn{13}{c|}{$S_{(12)}(i,m)$} \\
                 & $1$ & $2$ & $3$ & $4$ & $5$ & $6$ & $7$ & $8$ & $9$ & $10$ & $11$ & $12$ & $13$ \\
                 \hline
                1 & -$\frac{17}{108}$ & -$\frac{347}{108}$ & -$\frac{515}{108}$ & $\frac{43}{108}$ & $\frac{307}{108}$ & -$\frac{383}{108}$ & -$\frac{263}{108}$ & $\frac{511}{108}$ & $\frac{415}{108}$ & -$\frac{95}{108}$ & $\frac{211}{108}$ & $\frac{223}{108}$ & -$\frac{71}{108}$ \\ \hline 
                2 & -$\frac{127}{324}$ & $\frac{89}{324}$ & $\frac{845}{324}$ & $\frac{593}{324}$ & $\frac{593}{324}$ & $\frac{485}{324}$ & -$\frac{883}{324}$ & -$\frac{91}{324}$ & $\frac{557}{324}$ & -$\frac{523}{324}$ & -$\frac{1891}{324}$ & $\frac{89}{324}$ & $\frac{737}{324}$ \\ \hline 
                3 & -$\frac{95}{108}$ & -$\frac{71}{108}$ & -$\frac{515}{108}$ & -$\frac{383}{108}$ & $\frac{415}{108}$ & $\frac{223}{108}$ & -$\frac{347}{108}$ & $\frac{307}{108}$ & $\frac{511}{108}$ & $\frac{211}{108}$ & -$\frac{17}{108}$ & -$\frac{71}{108}$ & -$\frac{263}{108}$ \\ \hline 
                4 & $\frac{55}{108}$ & $\frac{55}{108}$ & $\frac{451}{108}$ & $\frac{451}{108}$ & -$\frac{263}{108}$ & -$\frac{155}{108}$ & $\frac{451}{108}$ & -$\frac{155}{108}$ & -$\frac{551}{108}$ & -$\frac{263}{108}$ & $\frac{55}{108}$ & -$\frac{263}{108}$ & -$\frac{155}{108}$ \\ \hline 
                5 & -$\frac{319}{324}$ & $\frac{617}{324}$ & -$\frac{1111}{324}$ & -$\frac{1687}{324}$ & $\frac{329}{324}$ & $\frac{509}{324}$ & -$\frac{175}{324}$ & $\frac{365}{324}$ & $\frac{545}{324}$ & $\frac{509}{324}$ & $\frac{959}{324}$ & $\frac{257}{324}$ & -$\frac{1039}{324}$  \\ \hline 
                6 & $\frac{1039}{324}$ & -$\frac{545}{324}$ & $\frac{175}{324}$ & $\frac{1543}{324}$ & -$\frac{257}{324}$ & -$\frac{617}{324}$ & -$\frac{329}{324}$ & -$\frac{617}{324}$ & -$\frac{437}{324}$ & -$\frac{257}{324}$ & -$\frac{707}{324}$ & -$\frac{509}{324}$ & $\frac{1507}{324}$  \\ \hline 
                7 & $\frac{257}{324}$ & $\frac{365}{324}$ & -$\frac{1687}{324}$ & -$\frac{751}{324}$ & $\frac{545}{324}$ & $\frac{329}{324}$ & -$\frac{319}{324}$ & $\frac{509}{324}$ & $\frac{509}{324}$ & $\frac{617}{324}$ & $\frac{959}{324}$ & -$\frac{175}{324}$ & -$\frac{1111}{324}$  \\ \hline 
                8 & -$\frac{17}{108}$ & $\frac{307}{108}$ & $\frac{415}{108}$ & -$\frac{71}{108}$ & -$\frac{71}{108}$ & $\frac{511}{108}$ & $\frac{223}{108}$ & -$\frac{515}{108}$ & -$\frac{263}{108}$ & $\frac{211}{108}$ & -$\frac{347}{108}$ & -$\frac{383}{108}$ & -$\frac{95}{108}$  \\ \hline 
                9 & -$\frac{617}{324}$ & -$\frac{257}{324}$ & -$\frac{545}{324}$ & -$\frac{617}{324}$ & -$\frac{437}{324}$ & -$\frac{617}{324}$ & $\frac{175}{324}$ & $\frac{1507}{324}$ & -$\frac{257}{324}$ & -$\frac{329}{324}$ & $\frac{1543}{324}$ & $\frac{1039}{324}$ & -$\frac{707}{324}$  \\ \hline 
                10 & -$\frac{257}{324}$ & -$\frac{617}{324}$ & $\frac{1543}{324}$ & $\frac{1507}{324}$ & -$\frac{437}{324}$ & -$\frac{257}{324}$ & $\frac{1039}{324}$ & -$\frac{257}{324}$ & -$\frac{617}{324}$ & -$\frac{545}{324}$ & -$\frac{707}{324}$ & -$\frac{329}{324}$ & $\frac{175}{324}$  \\ \hline 
                11 & $\frac{845}{324}$ & -$\frac{883}{324}$ & -$\frac{1891}{324}$ & $\frac{89}{324}$ & $\frac{485}{324}$ & -$\frac{523}{324}$ & -$\frac{127}{324}$ & $\frac{593}{324}$ & $\frac{557}{324}$ & $\frac{737}{324}$ & $\frac{593}{324}$ & -$\frac{91}{324}$ & $\frac{89}{324}$ \\ \hline 
                12 & $\frac{737}{324}$ & $\frac{845}{324}$ & $\frac{485}{324}$ & $\frac{557}{324}$ & $\frac{89}{324}$ & $\frac{89}{324}$ & $\frac{593}{324}$ & -$\frac{91}{324}$ & -$\frac{1891}{324}$ & -$\frac{127}{324}$ & $\frac{593}{324}$ & -$\frac{883}{324}$ & -$\frac{523}{324}$ \\ \hline 
                13 & -$\frac{1687}{324}$ & -$\frac{175}{324}$ & $\frac{509}{324}$ & -$\frac{751}{324}$ & -$\frac{1111}{324}$ & $\frac{509}{324}$ & $\frac{545}{324}$ & $\frac{257}{324}$ & $\frac{617}{324}$ & $\frac{329}{324}$ & $\frac{365}{324}$ & $\frac{959}{324}$ & -$\frac{319}{324}$  \\ \hline 
                14 & $\frac{211}{108}$ & -$\frac{95}{108}$ & $\frac{415}{108}$ & $\frac{511}{108}$ & -$\frac{263}{108}$ & -$\frac{383}{108}$ & $\frac{307}{108}$ & -$\frac{71}{108}$ & -$\frac{515}{108}$ & -$\frac{347}{108}$ & -$\frac{17}{108}$ & -$\frac{71}{108}$ & $\frac{223}{108}$  \\ \hline 
                15 & $\frac{737}{324}$ & $\frac{89}{324}$ & -$\frac{1891}{324}$ & -$\frac{523}{324}$ & $\frac{557}{324}$ & -$\frac{91}{324}$ & -$\frac{883}{324}$ & $\frac{485}{324}$ & $\frac{593}{324}$ & $\frac{593}{324}$ & $\frac{845}{324}$ & $\frac{89}{324}$ & -$\frac{127}{324}$ \\ \hline 
                16 & -$\frac{17}{108}$ & $\frac{211}{108}$ & $\frac{511}{108}$ & $\frac{307}{108}$ & -$\frac{347}{108}$ & $\frac{223}{108}$ & $\frac{415}{108}$ & -$\frac{383}{108}$ & -$\frac{515}{108}$ & -$\frac{71}{108}$ & -$\frac{95}{108}$ & -$\frac{263}{108}$ & $\frac{43}{108}$ \\ \hline 
                17 & -$\frac{263}{108}$ & $\frac{55}{108}$ & -$\frac{263}{108}$ & -$\frac{551}{108}$ & -$\frac{155}{108}$ & $\frac{451}{108}$ & -$\frac{155}{108}$ & -$\frac{263}{108}$ & $\frac{451}{108}$ & $\frac{451}{108}$ & $\frac{55}{108}$ & $\frac{55}{108}$ & -$\frac{155}{108}$ \\ \hline 
                18 & $\frac{959}{324}$ & $\frac{509}{324}$ & $\frac{545}{324}$ & $\frac{365}{324}$ & -$\frac{175}{324}$ & $\frac{509}{324}$ & $\frac{329}{324}$ & -$\frac{1687}{324}$ & -$\frac{1111}{324}$ & $\frac{617}{324}$ & -$\frac{319}{324}$ & -$\frac{751}{324}$ & $\frac{257}{324}$ \\ \hline 
                19 & -$\frac{707}{324}$ & -$\frac{257}{324}$ & -$\frac{437}{324}$ & -$\frac{617}{324}$ & -$\frac{329}{324}$ & -$\frac{617}{324}$ & -$\frac{257}{324}$ & $\frac{1543}{324}$ & $\frac{175}{324}$ & -$\frac{545}{324}$ & $\frac{1039}{324}$ & $\frac{1507}{324}$ & -$\frac{257}{324}$ \\ \hline 
                20 & $\frac{959}{324}$ & $\frac{617}{324}$ & $\frac{509}{324}$ & $\frac{509}{324}$ & -$\frac{319}{324}$ & $\frac{329}{324}$ & $\frac{545}{324}$ & -$\frac{1039}{324}$ & -$\frac{1687}{324}$ & $\frac{365}{324}$ & $\frac{257}{324}$ & -$\frac{1111}{324}$ & -$\frac{175}{324}$ \\ \hline 
                21 & -$\frac{95}{108}$ & -$\frac{383}{108}$ & -$\frac{347}{108}$ & $\frac{211}{108}$ & -$\frac{263}{108}$ & -$\frac{515}{108}$ & $\frac{223}{108}$ & $\frac{511}{108}$ & $\frac{43}{108}$ & -$\frac{71}{108}$ & $\frac{415}{108}$ & $\frac{307}{108}$ & -$\frac{17}{108}$ \\ \hline 
                22 & $\frac{1543}{324}$ & -$\frac{329}{324}$ & -$\frac{257}{324}$ & $\frac{1507}{324}$ & $\frac{175}{324}$ & -$\frac{617}{324}$ & -$\frac{437}{324}$ & $\frac{13}{324}$ & -$\frac{545}{324}$ & -$\frac{257}{324}$ & -$\frac{617}{324}$ & -$\frac{707}{324}$ & $\frac{1039}{324}$ \\ \hline 
                23 & -$\frac{509}{324}$ & $\frac{175}{324}$ & -$\frac{329}{324}$ & -$\frac{707}{324}$ & -$\frac{545}{324}$ & -$\frac{617}{324}$ & -$\frac{257}{324}$ & $\frac{1039}{324}$ & -$\frac{257}{324}$ & -$\frac{437}{324}$ & $\frac{1507}{324}$ & $\frac{1543}{324}$ & -$\frac{617}{324}$ \\ \hline 
                24 & $\frac{89}{324}$ & -$\frac{91}{324}$ & $\frac{593}{324}$ & $\frac{737}{324}$ & $\frac{557}{324}$ & $\frac{593}{324}$ & -$\frac{127}{324}$ & -$\frac{523}{324}$ & $\frac{485}{324}$ & $\frac{89}{324}$ & -$\frac{1891}{324}$ & -$\frac{883}{324}$ & $\frac{845}{324}$ \\ \hline 
                25 & -$\frac{883}{324}$ & $\frac{593}{324}$ & -$\frac{127}{324}$ & -$\frac{1891}{324}$ & -$\frac{91}{324}$ & $\frac{593}{324}$ & $\frac{89}{324}$ & $\frac{89}{324}$ & $\frac{557}{324}$ & $\frac{485}{324}$ & $\frac{845}{324}$ & $\frac{737}{324}$ & -$\frac{523}{324}$ \\ \hline 
                26 & $\frac{959}{324}$ & $\frac{365}{324}$ & $\frac{329}{324}$ & $\frac{617}{324}$ & $\frac{257}{324}$ & $\frac{545}{324}$ & $\frac{509}{324}$ & -$\frac{1111}{324}$ & -$\frac{1039}{324}$ & $\frac{509}{324}$ & -$\frac{175}{324}$ & -$\frac{1687}{324}$ & -$\frac{319}{324}$ \\ \hline 
                $\mathbf{\Sigma}_{(12)}$ & $\frac{64}{9}$ & $\frac{17}{18}$ & -$\frac{127}{18}$ & $\frac{113}{18}$ & -$\frac{113}{18}$ & $\frac{25}{18}$ & 8 & $\frac{5}{2}$ & -$\frac{251}{18}$ & $\frac{53}{9}$ & $\frac{221}{18}$ & -$\frac{64}{9}$ & -$\frac{47}{6}$ \\
                \hline
            \end{tabular}
            \caption{Computation of $S_{(12)}(i,m) = (-1)^m \sum_{\substack{\tiny{(l_1,1_2,l_3)\in\mathbb{Z}_{\geq0}^{3}}\\1_1 + l_2+l_3=3}} \prod_{k=1}^3 \frac{B_{l_k}(\tilde{x}_{(12),k}(i, m))}{l_k!} \,\tr_{F/\Q}\bigg(\prod_{k=1}^3 f_{(12),k}^{l_k-1}\bigg)$}
            \label{table:sums2}
        \end{table}
        \newpage
        \noindent Using \texttt{SageMath}, we compute $w_K=6$, $[\oa_F^{\times}:\oa_F^{\times, +}]=8$, and $[\oa_F^{\times, +}:N_{K/F}\oa_K^{\times}] =1.$ Adding the contributions from each $S_\mathrm{id}(i,m)$ listed in Table \ref{table:sums1} and the contributions from each $S_{(12)}(i,m)$ listed in Table \ref{table:sums2}, we compute that the class number is given by
        \begin{align*}
        h_K&=\frac{1}{4}\left(\frac{199}{36}+\frac{155}{36}+\frac{64}{9}+\frac{17}{18}-\frac{127}{18}+\frac{113}{18}-\frac{113}{18}+\frac{25}{18}+8+\frac{5}{2}-\frac{251}{18}+\frac{53}{9}+\frac{221}{18}-\frac{64}{9}-\frac{47}{6}\right)=\frac{12}{4} = 3.&
        \end{align*}
    \end{example}

\section*{Acknowledgements}
The author is a participant in the 2023 UVA REU in Number Theory. She is grateful for the support of grants from Jane Street Capital, the National Science Foundation (DMS-2002265 and DMS- 2147273), the National Security Agency (H98230-23-1-0016), and the Templeton World Charity Foundation. She thanks Ken Ono for suggesting this problem and for his mentorship and support, Wei-Lun Tsai for his mentorship and support, and Eleanor McSpirit and Alejandro De Las Pe\~nas Casta\~no for their mentorship and guidance. She would also like to thank Samit Dasgupta and the reviewer for their helpful and constructive comments, which have improved the presentation and quality of this paper.

\bibliography{bibliography}
\end{document}